\documentclass{article}
\usepackage{arxiv}

\usepackage[utf8]{inputenc} 
\usepackage[T1]{fontenc}    

\usepackage{hyperref}       
\usepackage{url}            
\usepackage{booktabs}       
\usepackage{amsfonts}       
\usepackage{nicefrac}       
\usepackage{microtype}      
\usepackage{lipsum}		
\usepackage{graphicx}
\usepackage[numbers]{natbib}
\usepackage{doi}
\usepackage{graphicx} 
\usepackage{amsmath, amssymb, bm, xcolor,mathtools}
\usepackage{algorithm}
\usepackage{algorithmic}
\usepackage{amsthm}

\numberwithin{equation}{section}

\newcommand{\defeq}{\vcentcolon=}

\newtheorem{theorem}{Theorem}[section] 
\newtheorem{lemma}[theorem]{Lemma}     
\newtheorem{proposition}[theorem]{Proposition}

\newtheorem{definition}[theorem]{Definition}

\title{MAGPIE: Multilevel--Adaptive--Guided solver for ptychographic phase retrieval}

\author{ \href{https://orcid.org/0009-0001-9020-9152}{\includegraphics[scale=0.06]{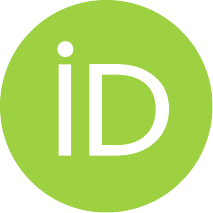}\hspace{1mm}Borong Zhang}\thanks{Corresponding author.} \\
	Department of Mathematics\\
	University of Wisconsin--Madison\\
	Madison, WI 53706 \\
	\texttt{bzhang388@wisc.edu} \\
	\And
	\href{https://orcid.org/0000-0001-9210-8948}{\includegraphics[scale=0.06]{orcid.pdf}\hspace{1mm}Qin Li} \\
	Department of Mathematics\\
	University of Wisconsin--Madison\\
	Madison, WI 53706 \\
	\texttt{qinli@math.wisc.edu} \\
    \And
	\href{https://orcid.org/0000-0002-4131-9363}{\includegraphics[scale=0.06]{orcid.pdf}\hspace{1mm}Zichao Wendy Di} \\
	Mathematics and Computer Science\\
	Argonne National Laboratory\\
	Lemont, IL 60439\\
	\texttt{wendydi@anl.gov} \\
}

\renewcommand{\shorttitle}{Stochastic Multigrid Ptychography}

\hypersetup{
pdftitle={MAGPIE: Multilevel--Adaptive--Guided solver for ptychographic phase retrieval},
pdfsubject={math, optics},
pdfauthor={Borong Zhang, Qin Li, Zichao Wendy Di},
pdfkeywords={phase retrieval,  multigrid optimization,  inverse problems,  ptychography, coherent diffraction imaging},
}

\begin{document}
\maketitle

\begin{abstract}
We introduce MAGPIE (Multilevel––Adaptive––Guided Ptychographic Iterative Engine), a stochastic multigrid solver for the ptychographic phase-retrieval problem. The ptychographic phase-retrieval problem is inherently nonconvex and ill-posed. To address these challenges, we reformulate the original nonlinear and nonconvex inverse problem as the iterative minimization of a quadratic surrogate model that majorizes the original objective. This surrogate not only ensures favorable convergence properties but also generalizes the Ptychographic Iterative Engine (PIE) family of algorithms. By solving the surrogate model using a multigrid method, MAGPIE achieves substantial gains in convergence speed and reconstruction quality over traditional approaches.
\end{abstract}

\keywords{phase retrieval \and multigrid optimization \and inverse problems \and ptychography \and coherent diffraction imaging}

\section{Introduction}
Ptychography is a state-of-the-art coherent diffractive imaging technique that has attracted significant attention for its ability to reconstruct high-resolution, complex-valued images from a series of diffraction patterns. The technique finds applications in materials science~\cite{crystallography, materials_science_2}, biological imaging~\cite{biology_1, biology_2}, and integrated circuit imaging~\cite{integrated_circuit, integrated_circuit_1}(see~\cite{Rodenburg2019} for a comprehensive overview). In a typical ptychographic experiment, a coherent beam is scanned over overlapping regions of an object while recording the resulting diffraction patterns. Since only the intensity of the diffracted waves is measured—and phase information is lost—the reconstruction task becomes a \emph{phase retrieval problem}, requiring the recovery of both the amplitude and the phase of the object.

Mathematically, ptychographic phase retrieval is formulated as an inverse problem that is inherently nonconvex and ill-posed. To address these challenges, numerous algorithms have been proposed, including the well-known Ptychographical Iterative Engine (PIE) family~\cite{PIE,ePIE,rPIE}. This family of algorithms casts the problem in an iterative fashion and stabilizes the reconstruction through adding a quadratic regularization vector. 

Another feature of the ptychographic phase retrieval problem is its inherent multiscale structure. Coarse and fine features coexist in the to-be-reconstructed object, requiring different levels of numerical accuracy. This naturally calls for a multigrid solver. Multigrid methods have a long history in the numerical solution of partial differential equations, exploiting the inherent multi-scale structure of the underlying problem by solving it on a hierarchy of discretizations~\cite{multigrid_book}. Coarse-grid approximations capture the global structure of the solution, while fine-grid corrections resolve detailed features. In the context of optimization, a multilevel\footnote{Throughout this work, the terms ``grid" and ``level" are used interchangeably.} optimization algorithm (termed MG/OPT)~\cite{MGOPT} has demonstrated significant improvements in convergence performance, while maintaining flexibility over a wide range of applications~\cite{Nash01012000, mgopt_1}. Specific to inverse problems, this hierarchical structure can also help mitigate the ill-posedness inherent in many inverse problems~\cite{mginverse1}. Recent work has adapted MG/OPT to the ptychographic context~\cite{MGOPT_ptycho}, resulting in reduced computational costs and enhanced convergence performance.

In this paper, we propose MAGPIE (Multilevel--Adaptive--Guided Ptychographic Iterative Engine), a solver that integrates the advantages of both PIE and MG/OPT. As done in PIE, we recast the problem into an iterative minimization problem with each step solving a quadratic surrogate that majorizes the original objective. The majorization effect brings us a favorable convergence rate. Meanwhile, we also deploy a multigrid-based optimization strategy and reconstruct the coarse and fine features of the object level by level. Unlike other multigrid-based optimization methods whose coarse-level models are generically incompatible with those on the fine level, our proposed surrogate naturally builds compatibility across scales. This compatibility also ensures automatic hyperparameter selection: once parameters are tuned on the finest grid, those for coarser levels are derived automatically, simplifying tuning in the implementation and improving usability. The overall method is run in a stochastic fashion, and the final algorithm delivers a nice reconstruction with substantially reduced costs of full-gradient evaluations and accelerated computational performance.

The remainder of this paper is organized as follows. Section~\ref{sec:prelim} presents the ptychographic phase retrieval problem and introduces two main techniques to be deployed: the surrogate modeling deployed in the PIE family algorithm and the multigrid optimization strategy; Section~\ref{sec:surrogate} introduces our surrogate model, a model that will be shown to majorize the original objective and thus brings in favorable convergence properties; in Section~\ref{sec:multigrid} we finally present our ultimate solver: MAGPIE, through designing a multigrid method for the surrogate model. We also examine the iterative performance of the algorithm. In Section~\ref{sec:numerical_examples} we present numerical results. The numerical implementation of {MAGPIE} is available online, together with reproducible experiments and instructions.\footnote{Repository URL: \url{https://github.com/borongzhang/magpie_ptychography}.}
\section{Preliminaries}\label{sec:prelim}

\subsection{2D ptychographic phase retrieval} 
\begin{figure}[htbp]
    \centering
    \includegraphics[width=0.8\textwidth]{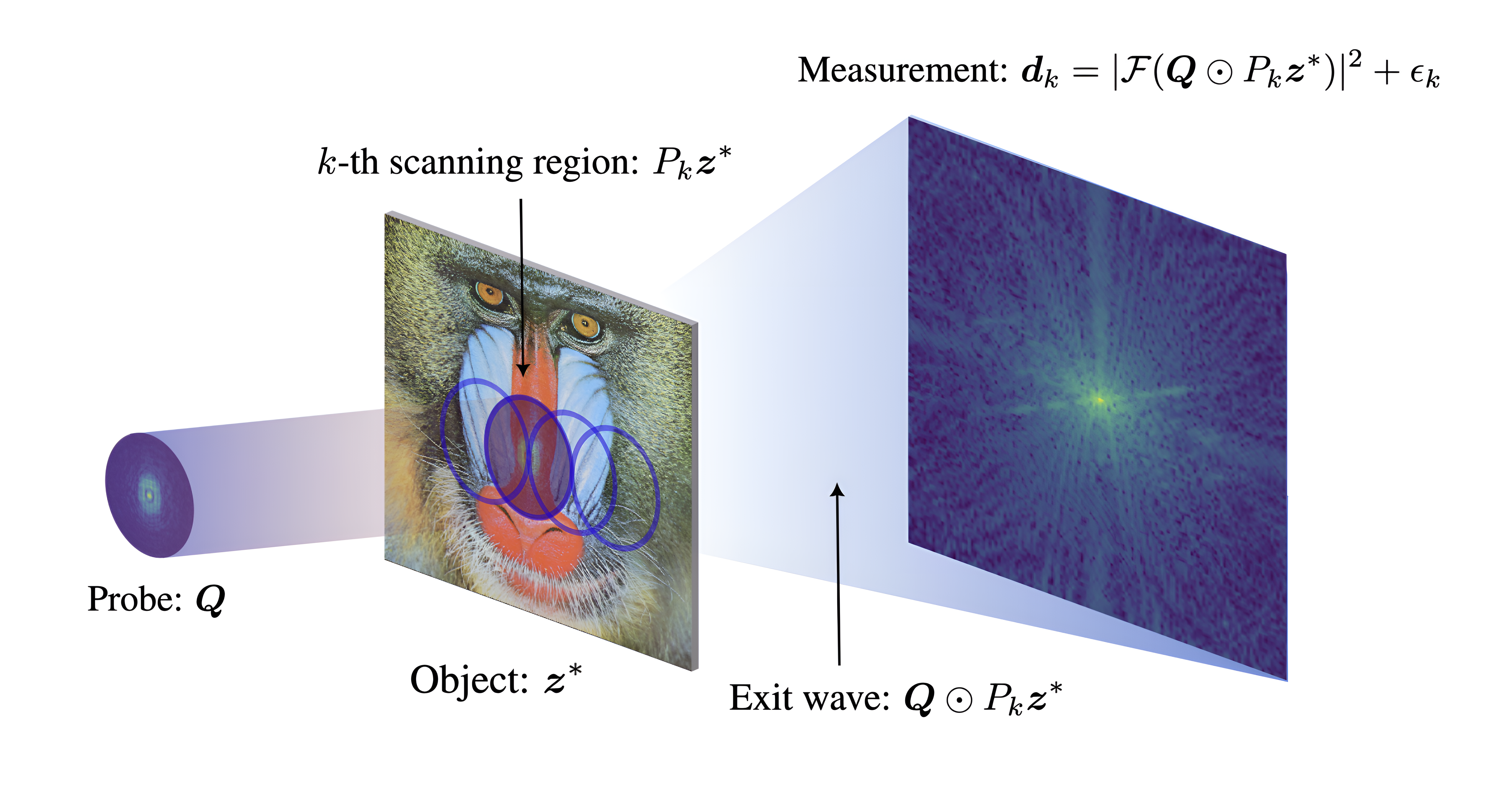}
    \caption{Experimental setup and data acquisition for ptychography.}
    \label{fig:setup}
\end{figure}

In 2D ptychographic phase retrieval, the goal is to reconstruct the object $\bm{z}^\ast \in \mathbb{C}^{n^2}$ (in vectorized form) from measured diffraction intensities. The data $\bm{d}_k \in \mathbb{R}^{m^2}$ are phaseless intensity measurements collected from overlapping probe positions:
\begin{equation*} \bm{d}_k = |\mathcal{F}(\bm{Q}\odot P_k \bm{z}^\ast)|^2 + \epsilon_k, \quad k = 1, \ldots, N, \end{equation*}
where $N$ is the number of scans, $P_k \in \mathbb{R}^{m^2\times n^2}$ extracts the $k$-th scanning region, $\bm{Q}\in\mathbb{C}^{m^2}$ is the probe, $\mathcal{F}$ is the two-dimensional discrete Fourier transform, $\epsilon_k$ models noise, $\odot$ denotes element-wise multiplication, and $|\cdot|^2$ applies element-wise to yield intensities. We refer to $\bm{Q}\odot P_k \bm{z}^\ast$ as the $k$-th exit wave. The overlap ratio is the fraction of the probe footprint shared between adjacent scans (e.g., $50\%$ overlap corresponds to half-probe shifts). Higher overlap ratios increase data redundancy and generally improve reconstruction stability and convergence, at the cost of additional measurements~\cite{Rodenburg2019}. Figure~\ref{fig:setup} illustrates the ptychographic acquisition geometry. It is worth noting that the probe $\bm{Q}$ can be highly nonuniform or concentrated, so some pixels in the $k$-th  scanning region receive much stronger illumination than others.

The inverse problem consists of inferring the object of interest, $\bm{z}^*$, from the collected data $\{\bm{d}_k\}_{k=1}^N$. This problem is commonly formulated as the following optimization problem:
\begin{equation*}
    \min_{\bm{z}} \Phi(\bm{z}) = \frac{1}{2} \sum_{k=1}^N \left\| |\mathcal{F}(\bm{Q}\odot P_k \bm{z})| - \sqrt{\bm{d}_k} \right\|_2^2,
\end{equation*}
where $\sqrt{\cdot}$ and $|\cdot|$ are applied element-wise. Following the terminology in~\cite{rPIE} and invoking Parseval's theorem, the problem can equivalently be written as:
\begin{equation}\label{eqn:exit_misfit}
\begin{aligned}
    \min_{\bm{z}}\Phi(\bm{z}) &= \frac{1}{2} \sum_{k=1}^N  \left\|\bm{Q}\odot P_k \bm{z} - \mathcal{R}_k(P_k \bm{z}) \right\|_2^2,\\
    &=\sum_{k=1}^N \Phi_k(\bm{z}_k) = \frac{1}{2} \sum_{k=1}^N \left\|\bm{Q}\odot\bm{z}_k - \mathcal{R}_k(\bm{z}_k) \right\|_2^2,
\end{aligned}
\end{equation}
where $\bm{z}_k$ and the operator $\mathcal{R}_k$ defining the revised exit wave are introduced below. {The equivalent formulation~\eqref{eqn:exit_misfit} rewrites the preceding amplitude-based data-misfit objective in terms of the exit wave variables used throughout the paper. This representation keeps the gradient formula, the PIE-type update, and the multigrid construction more compact and transparent, without changing the underlying optimization problem.}

\begin{definition}[Revised Exit Wave]\label{def:REW}
For the $k$-th scanning region, define
\begin{equation*}
\bm{z}_k = P_k \bm{z},
\end{equation*}
and define the revised exit wave as
\begin{equation*}
\mathcal{R}_k(\bm{z}_k) = \mathcal{F}^{-1}\left(\sqrt{\bm{d}_k} \odot \exp\left(i\theta\left(\mathcal{F}(\bm{Q}\odot \bm{z}_k)\right)\right)\right),
\end{equation*}
where $\theta(\cdot)$ denotes the element-wise principal argument.\footnote{Since $\theta(0)$ is undefined, at index $r$, we set $\left[\theta\left(\mathcal{F}(\bm{Q}\odot \bm{z}_k)\right)\right]_r = 0$ when $\left[\mathcal{F}(\bm{Q}\odot \bm{z}_k)\right]_r = 0$.}
\end{definition}

Since the objective function is real-valued but defined on a complex variable, it is not complex-differentiable in the classical sense~\cite{lang1985complex}. A common remedy is to employ the $\mathbb{C}\mathbb{R}$-calculus~\cite{CRcalculus}. To that end, we represent vectors in the complex field, $\bm{z}_k\in \mathbb{C}^{m^2}$, in the Euclidean space $\mathbb{R}^{2m^2}$ as follows:
\begin{equation*}
\left\{ [\bm{x}_k^\top, \bm{y}_k^\top]^\top : \bm{x}_k + i\bm{y}_k = \bm{z}_k \right\},
\end{equation*}
so that the objective function can be interpreted as a mapping from $\mathbb{R}^{2m^2}$ to $\mathbb{R}$. This representation enables us to compute the complex gradient of $\Phi_k$:
\begin{equation}\label{eqn:misfit_derivative}
\begin{aligned}
    \nabla_{\bm{z}_k} \Phi_k  &= \nabla_{\bm{x}_k}\Phi_k + i\nabla_{\bm{y}_k}\Phi_k,\\
    &= \overline{\bm{Q}}\odot\left(\bm{Q}\odot\bm{z}_k - \mathcal{R}_k(\bm{z}_k)\right),
\end{aligned}
\end{equation}
(see Appendix~\ref{app:computation_misfit_derivative} for the derivation). Consequently, the full gradient of $\Phi(\bm{z})$ in Eqn.~\eqref{eqn:exit_misfit} is
\begin{equation*} 
     \nabla_{\bm{z}}\Phi(\bm{z}) = \sum_{k=1}^N P_k^\top \nabla_{\bm{z}_k}\Phi_k(P_k\bm{z}) = \sum_{k=1}^N P_k^\top \left(\overline{\bm{Q}}\odot\left(\bm{Q}\odot P_k\bm{z} - \mathcal{R}_k\left(P_k\bm{z}\right)\right)\right).
\end{equation*}

\paragraph{PIE family of algorithms}
The PIE family of algorithms~\cite{PIE,ePIE,rPIE} is widely employed to minimize the misfit~ Eqn.~\eqref{eqn:exit_misfit}. At each iteration, the scanning regions are randomly shuffled and then updated one by one using a regularized quadratic surrogate model until all regions have been updated.  While one scanning region is selected and updated, all entries outside that region remain unchanged. Let $\bm{z}_k^{(j)}$ denote the estimate for the $k$-th region at the $j$-th iteration.\footnote{In each iteration, we shuffle the indices of the scanning regions $\{k\}_{k=1}^N$ and update them sequentially.} The local update is then performed by  
\begin{equation}\label{eqn:rpie_formulation}
\bm{z}_k^{(j)}\leftarrow \bm{z}_k^+ = \bm{z}_k^{(j)} + \overline{\bm{Q}}\oslash \left(\bm{u} + |\bm{Q}|^2\right)\odot\left(\mathcal{R}_k\left(\bm{z}_k^{(j)}\right) - \bm{Q}\odot\bm{z}_k^{(j)}\right),
\end{equation}
where $\oslash$ denotes the Hadamard (element-wise) division. Specifically, we  define its corresponding regularized quadratic surrogate function $\Phi^{\texttt{PIE}}_k\left(\bm{z}_k;\bm{z}_k^{(j)}\right)$ as
\begin{equation*} 
\Phi^{\texttt{PIE}}_k\left(\bm{z}_k ; \bm{z}_k^{(j)}\right) = \frac{1}{2} \left\|\bm{Q}\odot\bm{z}_k - \mathcal{R}_k\left(\bm{z}_k^{(j)}\right)\right\|_2^2 + \frac{1}{2}\bm{u}^\top\left|\bm{z}_k - \bm{z}_k^{(j)}\right|^2,
\end{equation*}
{where $\bm{u} \in \mathbb{R}^{m^2}$ is an a priori chosen regularization vector so that \(\bm{u}+|\bm{Q}|^2\) is strictly positive entrywise to ensure the Hadamard division is well-defined.} Then, $\bm{z}^{(j+1)}$ is obtained as the final output of  Eqn.~\eqref{eqn:rpie_formulation} after sweeping through all scanning regions. Notably, for a given $\bm{z}_k^{(j)}$, this objective function is quadratic in $\bm{z}_k$.
 
Correspondingly, the complex gradient of $\Phi^{\texttt{PIE}}_k\left(\bm{z}_k ; \bm{z}_k^{(j)}\right)$ with respect to $\bm{z}_k$ is given by
\begin{equation}\label{eqn:derivative_surrogate}
\nabla_{\bm{z}_k}\Phi^{\texttt{PIE}}_k\left(\bm{z}_k ; \bm{z}_k^{(j)}\right) = \overline{\bm{Q}}\odot\left(\bm{Q}\odot \bm{z}_k - \mathcal{R}_k\left(\bm{z}_k^{(j)}\right)\right) + \bm{u}\odot\left(\bm{z}_k - \bm{z}_k^{(j)}\right),
\end{equation}
(see Appendix~\ref{app:computation_surrogate_derivative} for the derivation). Therefore, setting the gradient to zero and solving for the explicit solution $\bm{z}_k^+$ recovers the local update rule in Eqn.~\eqref{eqn:rpie_formulation}.

The choice of $\bm{u}$ strongly influences the reconstruction accuracy. Within the PIE family of algorithms, experiments in~\cite{rPIE} indicate that the best performance is obtained by choosing
\begin{equation*}
    \bm{u} = \bm{u}^{\texttt{rPIE}}
    := \alpha\left(\|\bm{Q}\|^2_{\infty}\bm{J} - |\bm{Q}|^2\right),
\end{equation*}
where $\alpha > 0$, $\bm{J} \in \mathbb{R}^{m^2}$ denotes the vector with all entries equal to one, and $|\cdot|^2$ is applied element-wise. {For this choice, the denominator in Eqn.~\eqref{eqn:rpie_formulation} is strictly positive entrywise, and hence the Hadamard division is well defined.} This algorithm is termed the regularized Ptychographical Iterative Engine (rPIE). The rationale behind this choice is that at indices with low probe illumination ($|\bm{Q}_r|\ll 1$), the revised exit wave is more susceptible to noise, necessitating a larger penalty to stabilize the update. Conversely, where the probe strongly illuminates the object, a smaller penalty is applied, reflecting our increased confidence in those measurements. This adaptive weighting not only enhances the robustness of rPIE but also effectively mitigates noise-induced instability. The pseudocode for rPIE is presented in Algorithm~\ref{alg:rPIE}. However, despite these advantages, rPIE currently lacks a rigorous theoretical foundation in terms of convergence. {We reexamine the algorithm from a surrogate-minimization perspective in Section~\ref{sec:surrogate}.}

\begin{algorithm}
\caption{$\texttt{rPIE}(\bm{z}^{(0)},\bm{Q},\{\bm{d}_k\}_{k=1}^N,\alpha)$}\label{alg:rPIE}
\begin{algorithmic}[1]
\REQUIRE Object: $\bm{z}^{(0)}$, probe: $\bm{Q}$, data: $\{\bm{d}_k\}_{k=1}^N$, regularization constant: $\alpha$.
\STATE $\bm{u} \leftarrow \alpha\left(\|\bm{Q}\|^2_{\infty}\bm{J} - |\bm{Q}|^2\right)$
\WHILE{$\bm{z}^{(j)}$ not converged} 
    \FOR{{$k$ in a uniformly random permutation of $\{1,\dots,N\}$}}
        \STATE $\bm{z}_k^{(j)}\leftarrow \bm{z}_k^{(j)} + \overline{\bm{Q}}\oslash\left(\bm{u} + |\bm{Q}|^2\right)\odot\left(\mathcal{R}_k\left(\bm{z}_k^{(j)}\right) - \bm{Q}\odot\bm{z}_k^{(j)}\right)$
    \ENDFOR
\ENDWHILE
\RETURN $\bm{z}^{(j)}$
\end{algorithmic}
\end{algorithm}

\subsection{Multilevel optimization algorithm (MG/OPT)}\label{sec:MGOPT}
Ptychographic phase retrieval contains an intrinsic resolution hierarchy, making it an ideal candidate for multigrid acceleration. {MG/OPT is a multigrid solver for optimization. Given an objective function $f_h$ defined on a fine grid, our goal is to solve $\min_{\bm{z}} f_h(\bm{z})$. MG/OPT accelerates this process by, when appropriate, lifting the problem to a coarser level where a simpler optimization problem can be solved more efficiently. By communicating information between fine and coarse grids and offloading expensive fine-grid computations onto the cheaper coarse level, the method speeds up overall convergence. It has already demonstrated faster convergence for the ptychographic phase-retrieval problem~\cite{MGOPT_ptycho}. Here, we present the two-level MG/OPT algorithm, which employs one fine grid and one coarse grid; the multilevel extension follows naturally. We introduce the following notation:}
\begin{itemize}
  \item $h$ and $H$ denote the fine grid and coarse grid, respectively. Let $n_h$ be the number of fine-grid points per dimension. For simplicity, we assume $n_h = 2 n_H$.
  \item {$f_H$ is a manually defined objective function on the coarse grid, intended to capture key features of $f_h$.} 
    \item $\bm{I}_h^H$: The restriction operator (binning) that downsamples
    variables from the fine grid to the coarse grid via average pooling.
    Specifically, ${\bm{I}_h^H\in \mathbb{R}^{n_H^2 \times n_h^2}}$, and is given by
\begin{equation}\label{eqn:restriction_operator}
    \frac{1}{4}
    \begin{bmatrix}
    1 & 1 & 0 & \cdots & \cdots & \cdots & 0 & 1 & 1 & 0 & \cdots &  \cdots & \cdots & 0\\
    0 & 0 & 1 & 1 & 0 & \cdots & \cdots & \cdots & 0 & 1 & 1 & 0 & \cdots & 0\\
    &&&&\ddots&&&&\ddots &&&\ddots  \\
    0 & \cdots & \cdots & \cdots & 0 & 1 & 1 & 0 & \cdots & \cdots & \cdots & 0 & 1 & 1 
    \end{bmatrix}.
\end{equation}
  \item $\bm{I}_H^h$: The prolongation operator that interpolates the search direction from the coarse grid to the fine grid. Specifically, 
  \begin{equation}\label{eqn:prolongation_operator}
      \bm{I}_H^h = 4(\bm{I}_h^H)^\top.
  \end{equation}
\end{itemize}

We summarize the process of MG/OPT in Algorithm~\ref{alg:mgopt}.

\begin{algorithm}[H]
\caption{\texttt{MG/OPT}}
\label{alg:mgopt}
\begin{algorithmic}[1]
\REQUIRE Given an initial estimate of the solution $\bm{z}_h^0$ and integers $k_1, k_2 \ge 0$ satisfying $k_1 + k_2 > 0$.
\FOR{$j = 0, 1, \dots$ until convergence}
    \STATE \textbf{Pre-smoothing}: Apply $k_1$ iterations of a convergent optimization algorithm to $\min_{\bm{z}} f_h(\bm{z})$ to obtain $\bar{\bm{z}}_h$ (with $\bm{z}_h^{(j)}$ as the initial guess).
    \STATE \textbf{Recursion}:
    \begin{ALC@g}
        \STATE Compute $\bar{\bm{z}}_H = \bm{I}_h^H \bar{\bm{z}}_h$ and $\bar{\bm{v}}_H = \nabla f_H(\bar{\bm{z}}_H) - \bm{I}_h^H \nabla f_h(\bar{\bm{z}}_h)$.
        \STATE Minimize (perhaps approximately) the surrogate model
        \begin{equation*}
            f_s(\bm{z}_H) = f_H(\bm{z}_H) - \bar{\bm{v}}_H^\top \bm{z}_H
        \end{equation*}
        to obtain $\bm{z}_H^+$ (with $\bar{\bm{z}}_H$ as the initial guess). The minimization can be performed recursively by calling MG/OPT.
        \STATE Compute the search direction $\bm{e}_H = \bm{z}_H^+ - \bar{\bm{z}}_H$.
        \STATE Use a line search to determine $\bm{z}_h^+ = \bar{\bm{z}}_h + \alpha \bm{I}_H^h \bm{e}_H$ satisfying $f_h(\bm{z}_h^+) \le f_h(\bar{\bm{z}}_h)$.
    \end{ALC@g}
    \STATE \textbf{Post-smoothing}: Apply $k_2$ iterations of the same convergent optimization algorithm to $\min_{\bm{z}} f_h(\bm{z})$ to obtain $\bm{z}_h^{{(j+1)}}$ (with $\bm{z}_h^+$ as the initial guess).
\ENDFOR
\end{algorithmic}
\end{algorithm}

Note that the surrogate model optimized at the coarse-grid level is not the original $f_H$, but rather
\begin{equation*}
    f_s(\bm{z}_H) = f_H(\bm{z}_H)\underbrace{-\bar{\bm{v}}_H^\top\bm{z}_H}_{\text{correction term}},
\end{equation*}
where $-\bar{\bm{v}}_H^\top\bm{z}_H$, that we call the \emph{correction term}, enforces the matching of first-order conditions across levels:
\begin{equation*}
\nabla f_h(\bm{z}_h) = 0
\Longrightarrow
\bm{I}_h^H\nabla f_h(\bm{z}_h) = 0
\Longleftrightarrow
\nabla f_s\left(\bm{I}_h^H\bm{z}_h\right) = 0.
\end{equation*}

The rigorous descent analysis for MG/OPT is given in~\cite{MGOPT}. In the classical unconstrained setting, {under the hypotheses used in that analysis,}\footnote{{The hypotheses in such descent analyses include differentiability of $f_h$, compactness of the level set $\left\{\bm{x}_h:\ f_h(\bm{x}_h)\le f_h\left(\bm{x}_h^{(0)}\right)\right\}$, continuity of $\nabla^2 f_h$ on $S_h$, compatibility between the restriction and prolongation operators, and positive curvature of the coarse-grid objective $f_H$ along the coarse-grid correction direction $\bm{e}_H$.}} and with a sufficiently accurate coarse solve, the resulting fine-grid correction $\bm{e}_h$ is a descent direction for $f_h$ at $\bar{\bm{z}}_h$; that is,
\begin{equation*}
    \nabla f_h(\bar{\bm{z}}_h)^\top \bm{e}_h < 0.
\end{equation*}

Notably, the correction term ensures that the algorithm yields a descent direction even when $f_h$ and $f_H$ are inconsistent. However, the standard MG/OPT framework builds the coarse-grid surrogate $f_s$ by matching only the first-order (gradient) information of the fine-grid model. While this generic surrogate guarantees classical multigrid consistency, it neglects higher-order curvature and the strong block structure inherent in ptychography, causing convergence to slow down dramatically after only a few iterations. Instead, we propose a more robust coarse-level surrogate, consistent with the fine-level problem, which eliminates the need for a correction term and respects the ptychographic block structure. Section~\ref{sec:multigrid} presents this structure-aware surrogate and demonstrates its integration into our multigrid optimization framework.

\section{Surrogate Minimization}\label{sec:surrogate}
In this section, we introduce a surrogate objective for  Eqn.~\eqref{eqn:exit_misfit}. This surrogate not only majorizes  Eqn.~\eqref{eqn:exit_misfit} (see Section~\ref{sec:property_surrogate} for a precise definition of majorization) and iteratively recovers its optimal solution but is also sufficiently general to encompass the entire PIE family of solvers. This surrogate model has been studied in~\cite{MM_phaseretrieval} in a different context and has demonstrated superior performance on the general phase-retrieval problem.

Given $\bm{z}^{(j)}$ as the solution at the $j$-th iteration of an optimizer (e.g., rPIE), we define its corresponding quadratic surrogate model as:
\begin{equation}\label{eqn:quadratic_surrogate}
\begin{aligned}
    \widetilde{\Phi}\left(\bm{z} ; \bm{z}^{(j)}\right) &= \frac{1}{2} \sum_{k=1}^N \left\|\bm{Q}\odot P_k\bm{z} - \mathcal{R}_k\left(P_k\bm{z}^{(j)}\right)\right\|_2^2,\\
    &=\sum_{k=1}^N \widetilde{\Phi}_k\left(\bm{z}_k ; \bm{z}_k^{(j)}\right) = \frac{1}{2} \sum_{k=1}^N \left\|\bm{Q}\odot \bm{z}_k - \mathcal{R}_k\left(\bm{z}_k^{(j)}\right)\right\|_2^2,
\end{aligned}
\end{equation}
where $\mathcal{R}_k$ is defined in Definition~\ref{def:REW}. We note that in the definition of $\mathcal{R}_k$, the term $\left[\theta\left(\mathcal{F}\left(\bm{Q}\odot\bm{z}_k^{(j)}\right)\right)\right]_r$ is ill-defined if  $\left[\mathcal{F}\left(\bm{Q}\odot\bm{z}_k^{(j)}\right)\right]_r=0$ at some index $r$, so that we complement the original definition with the following specification:
\begin{equation}\label{eqn:specification}
\left[\theta\left(\mathcal{F}\left(\bm{Q}\odot\bm{z}_k^{(j)}\right)\right)\right]_r =
\begin{cases}
0, & \text{if } j=0 \& \left[\mathcal{F}\left(\bm{Q}\odot\bm{z}_k^{(j)}\right)\right]_r=0, \\
\left[\theta\left(\mathcal{F}\left(\bm{Q}\odot\bm{z}_k^{(j-1)}\right)\right)\right]_r, & \text{if } j\geq 1 \& \left[\mathcal{F}\left(\bm{Q}\odot\bm{z}_k^{(j)}\right)\right]_r=0,\\
\left[\theta\left(\mathcal{F}\left(\bm{Q}\odot\bm{z}_k^{(j)}\right)\right)\right]_r, & else.
\end{cases}
\end{equation}
In comparison to the original definition, when the new iteration does not provide any further information, instead of setting $\theta$, the principal argument directly to $0$, we recycle the value from the previous iteration, and thus is numerically more stable.

The $(j+1)$-th iteration solution is then a numerical solution to:
\begin{equation}\label{eqn:surrogate_minimization}
\bm{z}^{(j+1)} = \operatornamewithlimits{argmin}_{\bm{z}} \widetilde{\Phi}\left(\bm{z} ; \bm{z}^{(j)}\right).
\end{equation}

Moreover, observing that~ Eqn.~\eqref{eqn:quadratic_surrogate} is expressed as a summation, it is natural to adopt a stochastic optimization strategy. Specifically, in each iteration, we shuffle the indices of the scanning regions $\{k\}_{k=1}^N$ and update them sequentially:
\begin{equation}\label{eqn:stochastic_surrogate_minimization}
    \bm{z}_k^{(j)}\leftarrow\bm{z}_k^+=\operatornamewithlimits{argmin}_{\bm{z}_k}\widetilde{\Phi}_k\left(\bm{z}_k; \bm{z}_k^{(j)}\right).
\end{equation}

One desirable feature of the surrogate model~ Eqn.~\eqref{eqn:quadratic_surrogate} is that it is quadratic. Consequently, an explicit analytical solution is available. To avoid singularity caused by small or even near zero entries of $\bm{Q}$, we further relax $\widetilde{\Phi}_k$ by adding a regularization vector (similar to the PIE approach):
\begin{equation}\label{eqn:proximal_point_method}
\bm{z}_k^{(j)}\leftarrow\bm{z}_k^+ = \operatornamewithlimits{argmin}_{\bm{z}_k}\widetilde{\Phi}_k\left(\bm{z}_k; \bm{z}_k^{(j)}\right)+\frac{1}{2}\bm{u}^\top\left|\bm{z}_k-\bm{z}_k^{(j)}\right|^2.
\end{equation}

Notice that the update formulation~ Eqn.~\eqref{eqn:proximal_point_method} is equivalent to applying one step of the proximal point method~\cite{Proximal} to~ Eqn.~\eqref{eqn:stochastic_surrogate_minimization}.

\subsection{Properties of the surrogate model}\label{sec:property_surrogate}
\begin{figure}[htbp]
    \centering
    \includegraphics[width=0.5\textwidth]{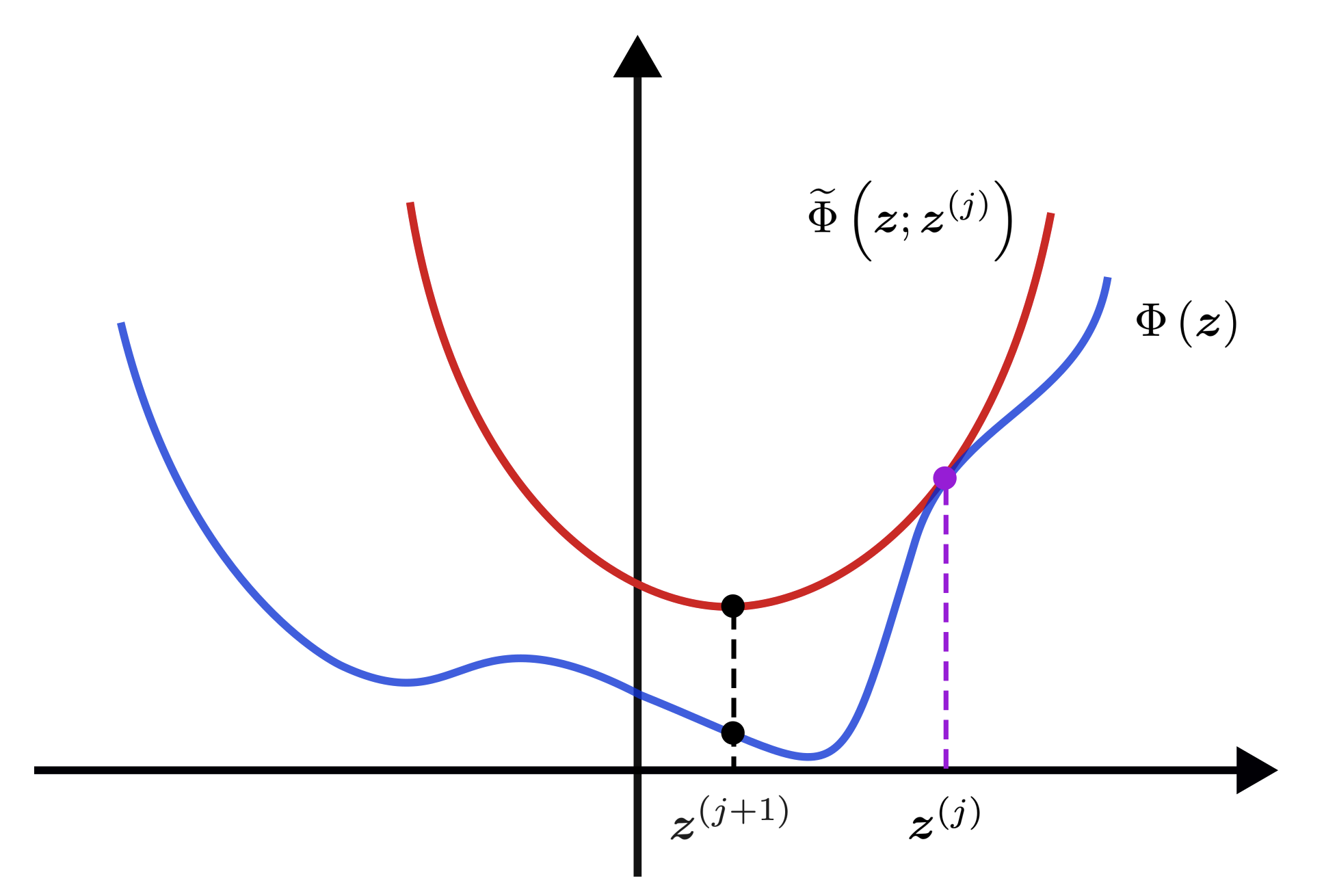}
    \caption{Illustration of the majorization property adapted from~\cite{SMM_graph}.}
    \label{fig:majorization}
\end{figure}

Our proposed surrogate model in  Eqn.~\eqref{eqn:surrogate_minimization} enjoys many desirable properties and is designed in accordance with the general principle of majorization~\cite{mm_alg_Lange} (see Figure~\ref{fig:majorization}).

\begin{proposition}[Majorization]\label{prop:mm_cond}
    The quadratic surrogate $\widetilde{\Phi}\left(\bm{z} ; \bm{z}^{(j)}\right)$, as an approximation to the original objective function $\Phi(\bm{z})$ in~ Eqn.~\eqref{eqn:exit_misfit}, has the following properties:
\begin{itemize}
    \item $\Phi(\bm{z})$ and $\widetilde{\Phi}\left(\bm{z} ; \bm{z}^{(j)}\right)$ agree at $\bm{z}^{(j)}$:
    \begin{equation*}
        \Phi\left(\bm{z}^{(j)}\right) = \widetilde{\Phi}\left(\bm{z}^{(j)} ; \bm{z}^{(j)}\right).
    \end{equation*}
    \item $\Phi(\bm{z})$ is dominated by $\widetilde{\Phi}\left(\bm{z} ; \bm{z}^{(j)}\right)$:
    \begin{equation}\label{eqn:major}
        \Phi\left(\bm{z}\right) \leq \widetilde{\Phi}\left(\bm{z} ; \bm{z}^{(j)}\right)\quad\text{for all}\quad\bm{z}\in\mathbb{C}^{n^2}.
    \end{equation}
    \item The complex gradients of $\Phi(\bm{z})$ and $\widetilde{\Phi}\left(\bm{z} ; \bm{z}^{(j)}\right)$ agree at $\bm{z}^{(j)}$:
    \begin{equation*} 
        \nabla_{\bm{z}}\Phi\left(\bm{z}^{(j)}\right) = \nabla_{\bm{z}}\widetilde{\Phi}\left(\bm{z}^{(j)} ; \bm{z}^{(j)}\right).
    \end{equation*}
\end{itemize}
\end{proposition}
The detailed proof of this proposition appears in Appendix~\ref{app:proof_mm_cond}. The first statement is straightforward to verify. To establish inequality~\eqref{eqn:major}, we use the definitions in Eqs.~\eqref{eqn:exit_misfit} and~\eqref{eqn:quadratic_surrogate} and verify the result for each $\Phi_k$. A similar proof can be found in~\cite{MM_phaseretrieval}. The proof of the final statement is adapted from~\cite{mm_convergence}.

Furthermore, minimizing the surrogate (see~ Eqn.~\eqref{eqn:surrogate_minimization}) provides a descent direction for the original objective function, as implied in the following proposition:
\begin{proposition}[{Monotone descent and stationarity}]\label{thm:convergence_iterative_surrogate}
Suppose that the quadratic surrogate $\widetilde{\Phi}\left(\bm{z} ; \bm{z}^{(j)}\right)$, given in  Eqn.~\eqref{eqn:quadratic_surrogate} has a global minimum $\bm{z}^{(j+1)}$. Then, for all $j=0,1,\dots$, we have
\begin{equation*}
    \Phi\left(\bm{z}^{(j+1)}\right)\leq\Phi\left(\bm{z}^{(j)}\right).
\end{equation*}
Additionally,
\begin{equation*}
    \Phi\left(\bm{z}^{(j+1)}\right)=\Phi\left(\bm{z}^{(j)}\right)\Longrightarrow \nabla_{\bm{z}}\Phi\left(\bm{z}^{(j+1)}\right) = 0.
\end{equation*}
\end{proposition}
The detailed proof of this proposition appears in Appendix~\ref{app:proof_convergence_iterative_surrogate}. The monotonicity inequality follows directly from Proposition~\ref{prop:mm_cond}, and the first-order optimality of $\Phi$ follows from the corresponding first-order optimality condition on the quadratic surrogate, with $\mathcal{R}_k$ specified in~\eqref{eqn:specification}. This is a standard result of a majorization–minimization algorithm~\cite{mm_convergence}. 

The next proposition gives a stationarity estimate for the surrogate
minimization in Eqn.~\eqref{eqn:surrogate_minimization}.

{\begin{proposition}[Sublinear stationarity bound]\label{thm:convergence_rate}
Let $L\defeq \frac12\left\|\sum_{k=1}^N P_k^\top |\bm{Q}|^2\right\|_\infty$.
Then,
\[
\frac{1}{4L}\sum_{j=0}^n
\left\|\nabla_{\bm z}\Phi\left(\bm z^{(j)}\right)\right\|_2^2
\le
\Phi\left(\bm{z}^{(0)}\right).
\]
Hence,
\[
\left\|\nabla_{\bm{z}}\Phi\left(\bm{z}^{(j)}\right)\right\|_2 \to 0,
\]
and in particular,
\[
\min_{0\le j\le n}
\left\|\nabla_{\bm{z}}\Phi\left(\bm{z}^{(j)}\right)\right\|_2^2
\le
\frac{4L}{n+1}\Phi\left(\bm{z}^{(0)}\right).
\]
\end{proposition}}
The detailed proof of this proposition appears in Appendix~\ref{app:proof_convergence_rate}. The overall derivation follows that in~\cite{mm_convergence}, where the increments at each iteration are bounded, and a telescoping-sum argument is applied.

\section{MAGPIE}\label{sec:multigrid}
As demonstrated in~\cite{MGOPT_ptycho}, the multigrid solver improves computational efficiency, while the surrogate majorization model provides favorable convergence properties. We now integrate these approaches into a single, tractable, and robust method, which we term the \emph{Multilevel–Adaptive–Guided Ptychographic Iterative Engine} (MAGPIE). The name reflects its multigrid structure, adaptive construction of coarse-grid surrogates, and guided fine‐grid reconstruction via coarse‐grid search directions. We present the algorithm in Section~\ref{sec:MAGPIE} and defer discussion of its well-definedness, consistency, descent properties, and automatic regularization selection to Section~\ref{sec:weights_property}.

\subsection{Algorithmic Development of MAGPIE}\label{sec:MAGPIE}
In this section, we present the proposed algorithm in a two-level formulation, which can be readily extended to multiple levels. Since there are only two levels, we omit the subscript $h$ from all terms defined on fine-grids.

For the $k$-th scanning region at the $j$-th iteration, the fine-grid surrogate model is:
\begin{equation}\label{eqn:fine_grid_surrogate_model}
\min_{\bm{z}_k}\widetilde{\Phi}_k\left(\bm{z}_k;\bm{z}_k^{(j)}\right) = \frac{1}{2}\left\|\bm{Q}\odot\bm{z}_k - \mathcal{R}_k\left(\bm{z}_k^{(j)}\right)\right\|_2^2,
\end{equation}
and we define the associated coarse-grid surrogate model as:
\begin{equation}\label{eqn:coarse_grid_surrogate_model}
\min_{\bm{z}_{H,k}} \widetilde{\Phi}_{H,k}\left(\bm{z}_{H,k};\bm{z}_k^{(j)}\right) = \frac{1}{2}\left\|\bm{Q}_H\odot\bm{z}_{H,k} - \mathcal{R}_{H,k}\left(\bm{z}_k^{(j)}\right)\right\|_2^2.
\end{equation}

The coarse-grid model represents a weighted average of the fine-grid model, with precise definitions given below:
\begin{itemize}
    \item The restriction operator $\bm{I}_h^H$ and the prolongation operator $\bm{I}_H^h$ are defined in  Eqn.~\eqref{eqn:restriction_operator} and  Eqn.~\eqref{eqn:prolongation_operator} respectively.  
    \item $\bm{Q}_H = \bm{I}_h^H\bm{Q}\in\mathbb{C}^{(m/2)^2}$: downsampled probe.
    \item $\bm{W}_{\mathcal{R}}\in\mathbb{C}^{m^2}$: {constant} weight vector for downsampling $\mathcal{R}_k\left(\bm{z}_k^{(j)}\right)$.
    \item $\mathcal{R}_{H,k}\left(\bm{z}_k^{(j)}\right) = \bm{I}_h^H\left(\bm{W}_{\mathcal{R}}\odot\mathcal{R}_k\left(\bm{z}_k^{(j)}\right)\right)$: downsampled revised exit wave.
\end{itemize}

Both the fine- and coarse-grid subproblems are quadratic and, in principle, admit closed-form solutions.  In practice, however, a direct solve is ill-posed: whenever entries of $\bm{Q}$ or $\bm{Q}_H$ vanish, the corresponding normal equations lose rank, and the solution is no longer unique.  As with the full surrogate problem in  Eqn.~\eqref{eqn:stochastic_surrogate_minimization}, we therefore replace an exact solve with a single proximal-point (relaxation) step. This involves the last iteration object estimate $\bm{z}_k^{(j)}$, its downsampled version $\bm{z}_{H,k}^{(j)}$, and two regularization vectors $\bm{u}$ (on fine-grid) and $\bm{u}_H$ (on coarse-grid). Concretely, we define:
\begin{itemize}
    \item $\bm{W}_{\bm{z}}\in \mathbb{R}^{m^2}_{\geq 0}$: {constant} weight vector for downsampling $\bm{z}_k^{(j)}$.
    \item {$\bm{z}_{H,k}^{(j)} = \bm{I}_h^H\left(\bm{W}_{\bm{z}}\odot\bm{z}_k^{(j)}\right)$: downsampled current object estimate.}
    \item $\bm{W}_{\bm{u}}^H\in \mathbb{R}^{(m/2)^2}_{\geq 0}$: {constant} weight vector for downsampling $\bm{u}$. 
    \item $\bm{u}_H = \bm{W}_{\bm{u}}^H \odot \bm{I}_h^H\bm{u}$: the coarse-grid regularization vector.
\end{itemize}

With these notations in place, we can proceed with updates. This goes through two steps. First, we derive the update for the coarse-grid model using the proximal point method in  Eqn.~\eqref{eqn:coarse_grid_surrogate_model}:
\begin{equation*} 
\hat{\bm{z}}_{H,k}^{(j)} = \bm{z}_{H,k}^{(j)} + \overline{\bm{Q}_H}\oslash\left(\bm{u}_H + |\bm{Q}_H|^2\right)\odot\left(\mathcal{R}_{H,k}\left(\bm{z}_k^{(j)}\right) - \bm{Q}_H\odot \bm{z}_{H,k}^{(j)}\right).
\end{equation*}
Correspondingly, the fine-grid update drawn from the coarse level is:
\begin{equation}\label{eqn:coarse-grid update}
    \widetilde{\bm{z}}_k^{(j)} = \bm{z}_k^{(j)} + \bm{I}_H^h\underbrace{\left(\hat{\bm{z}}_{H,k}^{(j)} - \bm{z}_{H,k}^{(j)}\right)}_{\bm{e}^{(j)}_{H,k}} = \bm{z}_k^{(j)} + \bm{I}_H^h \bm{e}^{(j)}_{H,k}.
\end{equation}

It is then followed by the second step where we apply the proximal point method on~ Eqn.~\eqref{eqn:fine_grid_surrogate_model} at the fine level to obtain:
\begin{equation*} 
\bm{z}_k^{(j)} \leftarrow \bm{z}_k^+ = \widetilde{\bm{z}}_k^{(j)} + \overline{\bm{Q}}\oslash\left(\bm{u} + |\bm{Q}|^2\right)\odot\left(\mathcal{R}_k\left(\bm{z}_k^{(j)}\right) - \bm{Q}\odot\widetilde{\bm{z}}_k^{(j)}\right).
\end{equation*}
In practice, we have the flexibility to choose the weights ($\bm{W}_{\bm{z}}$, $\bm{W}_{\mathcal{R}}$, and $\bm{W}_{\bm{u}}^H$). They are the coefficients that transfer information across levels. Therefore, they must be chosen carefully to ensure compatibility. We construct these weights in a sequential manner to guarantee consistency, well-definedness, descent directions, and stability.

\paragraph{Consistency}
Let the object estimate $\bm{\zeta}_k$ denote a minimizer of the fine-grid problem in Eqn.~\eqref{eqn:fine_grid_surrogate_model}. By consistency we mean that the downsampled object estimate $\bm{\zeta}_{H,k} = \bm{I}_h^H\left(\bm{W}_{\bm{z}}\odot\bm{\zeta}_k\right)$ should be a minimizer of the coarse-grid problem in Eqn.~\eqref{eqn:coarse_grid_surrogate_model}. We find that, in order to achieve this, the weight for downsampling the revised exit wave must satisfy
\[
\bm{W}_{\mathcal{R}} = (\bm{I}_H^h\bm{Q}_H)\odot\bm{W}_{\bm{z}}\oslash\bm{Q}.
\]
This relation is proved in Proposition~\ref{prop:consistency}.

\paragraph{Well-definedness}
The consistency relation above implies that the expression for $\bm{W}_{\mathcal{R}}$ may suffer from singularities due to division by $\bm{Q}$. To ensure well-definedness, we therefore choose $\bm{W}_{\bm{z}}$ so that these singularities are canceled. There are many possible choices; in this paper we take
\[
\bm{W}_{\bm{z}} = |\bm{Q}|^2\oslash\left(\bm{I}_H^h\bm{I}_h^H|\bm{Q}|^2\right),
\]
which has a desirable physical interpretation that we discuss in detail in Section~\ref{sec:weights_property}. In addition, Proposition~\ref{prop:weights_property} shows that all three weights are well defined and introduce no singularities.

\paragraph{Descent direction and stability}
We next need to show that the coarse-grid update in Eqn.~\eqref{eqn:coarse-grid update} gives a descent direction. In Theorem~\ref{prop:descent_direction}, we prove that this update yields a descent direction. We then control the magnitude of the update by tuning the coarse‐grid regularization vector. The regularization must balance stability (which favors a larger value) against the step size (which favors a smaller value). In Proposition~\ref{prop:automatic_regularization_selection}, we show that
\[
\bm{W}_{\bm{u}}^H = |\bm{Q}_H|^2\oslash\left(\bm{I}_h^H|\bm{Q}|^2\right)
\]
is a suitable choice.

We term this multigrid solver for each surrogate model the \emph{Multilevel--Adaptive--Guided Proximal Solver} (MAGPS). Specifically, ``Adaptive'' denotes the automatic construction of coarse‐grid surrogate terms using the weights specified in  Eqn.~\eqref{eqn:def_weights} without manual tuning, ``Guided'' emphasizes the directional coarse‐to‐fine corrections in  Eqn.~\eqref{eqn:coarse-grid update}, and ``Proximal Solver'' highlights the use of proximal‐point updates within the multigrid framework. The pseudocode for this solver is summarized in Algorithm~\ref{alg:MAGPS}.

The overall approach that addresses ptychographic phase retrieval by iteratively calling MAGPS to minimize the surrogate model over randomly sampled scanning regions is termed \emph{Multilevel--Adaptive--Guided Ptychographic Iterative Engine} (MAGPIE); its pseudocode is presented in Algorithm~\ref{alg:MAGPIE}. 

\begin{algorithm}
\caption{$\texttt{MAGPS}\left(\bm{z}^{(j)}_k, \bm{Q}, \mathcal{R}_k\left(\bm{z}_k^{(j)}\right),\bm{u}\right)$}
\label{alg:MAGPS}
\begin{algorithmic}[1]
\REQUIRE  Object $\bm{z}^{(j)}_k$, probe $\bm{Q}$, revised exit wave $\mathcal{R}_k\left(\bm{z}_k^{(j)}\right)$, regularization vector $\bm{u}$.
\IF{current level is coarsest}
    \STATE 
    $\bm{z}^+_k \leftarrow \bm{z}_k^{(j)} + \overline{\bm{Q}}\oslash\left(\bm{u} + |\bm{Q}|^2\right)\odot\left(\mathcal{R}_k\left(\bm{z}_k^{(j)}\right) - \bm{Q}\odot\bm{z}_k^{(j)}\right)$.
\ELSE
    \STATE \textbf{Coarse-Level Preparation:}
    \STATE Coarse probe: $\bm{Q}_H \leftarrow \bm{I}_h^H\bm{Q}$.
    \STATE {Set weights for downsampling the object:} $\bm{W}_{\bm{z}} \leftarrow |\bm{Q}|^2\oslash\left(\bm{I}_H^h\bm{I}_h^H|\bm{Q}|^2\right)$.
    \STATE {Set weights for downsampling the revised exit wave:} $\bm{W}_{\mathcal{R}} \leftarrow (\bm{I}_H^h\bm{Q}_H)\odot\bm{W}_{\bm{z}}\oslash\bm{Q}$.
    \STATE {Set weights for coarse-grid regularization:} $\bm{W}_{\bm{u}}^H \leftarrow |\bm{Q}_H|^2\oslash\left(\bm{I}_h^H|\bm{Q}|^2\right)$.
    \STATE Coarse object: $\bm{z}^{(j)}_{H,k} \leftarrow \bm{I}_h^H\left(\bm{W}_{\bm{z}}\odot\bm{z}^{(j)}_k\right)$.
    \STATE Coarse revised exit waves: $\mathcal{R}_{H,k}\left(\bm{z}_k^{(j)}\right) \leftarrow \bm{I}_h^H\left(\bm{W}_{\mathcal{R}}\odot\mathcal{R}_k\left(\bm{z}_k^{(j)}\right)\right)$.
    \STATE Coarse-grid regularization: $\bm{u}_H \leftarrow \bm{W}_{\bm{u}}^H\odot\bm{I}_h^H\bm{u}$.
    \STATE \textbf{Recursive Coarse Solve:}
    \STATE $\widehat{\bm{z}}^{(j)}_{H,k} \leftarrow \texttt{MAGPS}\left(\bm{z}^{(j)}_{H,k}, \bm{Q}_H,  \mathcal{R}_{H,k}\left(\bm{z}_k^{(j)}\right), \bm{u}_H\right)$.
    \STATE \textbf{Coarse-to-Fine Update:}
    \STATE Coarse-grid update direction: $\bm{e}^{(j)}_{H,k} \leftarrow \widehat{\bm{z}}^{(j)}_{H,k} - \bm{z}^{(j)}_{H,k}$.
    \STATE Update fine-grid estimate: $\widetilde{\bm{z}}^{(j)}_{k} \leftarrow \bm{z}^{(j)}_k + \bm{I}_H^h\bm{e}^{(j)}_{H,k}$.
    \STATE \textbf{Post-Smoothing:}
    \STATE $\bm{z}^+_k \leftarrow \widetilde{\bm{z}}^{(j)}_{k} + \overline{\bm{Q}}\oslash\left(\bm{u} + |\bm{Q}|^2\right)\odot\left(\mathcal{R}_k\left(\bm{z}_k^{(j)}\right) - \bm{Q}\odot\widetilde{\bm{z}}^{(j)}_{k}\right)$.
\ENDIF
\RETURN $\bm{z}^+_k$
\end{algorithmic}
\end{algorithm}

\begin{algorithm}
\caption{$\texttt{MAGPIE}(\bm{z}^{(0)},\bm{Q},\{\bm{d}_k\}_{k=1}^N,\bm{u})$}\label{alg:MAGPIE}
\begin{algorithmic}[1]
\REQUIRE Object: $\bm{z}^{(0)}$, probe: $\bm{Q}$, data: $\{\bm{d}_k\}_{k=1}^N$, regularization vector: $\bm{u}$.
\WHILE{$\bm{z}^{(j)}$ not converged} 
    \FOR{{$k$ in a uniformly random permutation of $\{1,\dots,N\}$}}
        \STATE $\bm{z}_k^{(j)}\leftarrow \texttt{MAGPS}_k\left(\bm{z}^{(j)}_k, \bm{Q}, \mathcal{R}_k\left(\bm{z}^{(j)}_k\right),\bm{u}\right)$
    \ENDFOR
\ENDWHILE
\RETURN $\bm{z}^{(j)}$
\end{algorithmic}
\end{algorithm}

\subsection{Converging properties of MAGPIE and weight selection}\label{sec:weights_property}

We summarize our choices of the weight vectors below. These definitions of $\bm{W}_{\bm{z}}$, $\bm{W}_{\mathcal{R}}$, and $\bm{W}_{\bm{u}}^H$ ensure the robust performance of MAGPIE.

\begin{equation}\label{eqn:def_weights}
    \bm{W}_{\bm{z}} = |\bm{Q}|^2\oslash\left(\bm{I}_H^h\bm{I}_h^H|\bm{Q}|^2\right),\quad
    \bm{W}_{\mathcal{R}} = (\bm{I}_H^h\bm{Q}_H)\odot\bm{W}_{\bm{z}}\oslash\bm{Q},\quad
    \bm{W}_{\bm{u}}^H = |\bm{Q}_H|^2\oslash\left(\bm{I}_h^H|\bm{Q}|^2\right).
\end{equation}

As shown in Propositions~\ref{prop:weights_property}, \ref{prop:consistency}, and~\ref{prop:automatic_regularization_selection}, and Theorem~\ref{prop:descent_direction}, these specific forms not only prevent singularities, maintain coarse-to-fine consistency, and preserve the descent property, but also ensure stability. It is worth noting that the quantity $\bm{I}_H^h\bm{I}_h^H|\bm{Q}|^2$
appears in both $\bm{W}_{\bm{z}}$ and $\bm{W}_{\mathcal{R}}$.
Its physical meaning is straightforward: within each $2\times 2$ bin of
$|\bm{Q}|^2$, the restriction $\bm{I}_h^H$ computes the average of the
four entries, and the prolongation $\bm{I}_H^h$ assigns this average back
to the original four fine-grid positions. By contrast, $|\bm{Q}|^2$
retains the original entrywise intensities. As a result,
\begin{equation*}
    \bm{z}_{H,k}^{(j)} = \bm{I}_h^H\left(\bm{W}_{\bm{z}}\odot\bm{z}_k^{(j)}\right)
\end{equation*} 
performs a bin-wise weighted average of $\bm{z}_k^{(j)}$, assigning greater weight to regions where $\bm{Q}$ has higher magnitude (i.e., higher signal-to-noise ratio). 

We first show that all weight vectors are element-wise bounded. 
\begin{proposition}[Well-definedness]\label{prop:weights_property}
{The element-wise magnitudes of the weight vectors defined in~Eqn.\eqref{eqn:def_weights} are bounded.} In particular:
\begin{equation*}
\left\|\bm{W}_{\bm{z}}\right\|_{\infty} \leq 4,\quad \left\|\bm{W}_{\mathcal{R}}\right\|_{\infty} \leq 4,\quad\text{and}\quad \left\|\bm{W}_{\bm{u}}^H\right\|_{\infty} \leq 1,
\end{equation*}
{where each weight vector's infinity norm equals the maximum magnitude among its entries.}
\end{proposition}

We leave the proof of this proposition to Appendix~\ref{app:proof_weights_property}.

We now show that the coarse-grid surrogate model is consistent with the fine-grid surrogate model. 
\begin{proposition}[Consistency]\label{prop:consistency}
Given the weight vectors defined in~ Eqn.~\eqref{eqn:def_weights}, the coarse-level objective $\widetilde{\Phi}_{H,k}$ (as defined in~ Eqn.~\eqref{eqn:coarse_grid_surrogate_model}) is consistent with the fine-level objective $\widetilde{\Phi}_k$ (as defined in~ Eqn.~\eqref{eqn:fine_grid_surrogate_model}) in the sense that, for all $\bm{z}_k\in\mathbb{C}^{m^2}$, given $\bm{z}_{H,k}=\bm{I}_h^H\left(\bm{W}_{\bm{z}}\odot\bm{z}_{k}\right)$, we have:
\begin{itemize}
    \item $\widetilde{\Phi}_{H,k}\left(\bm{z}_{H,k}; \bm{z}_k^{(j)}\right)\leq \frac{1}{4}\left\|\bm{W}_{\mathcal{R}}\right\|_{\infty}^2\widetilde{\Phi}_k\left(\bm{z}_k; \bm{z}_k^{(j)}\right)$.
    \item $\left\|\nabla_{\bm{z}_{H,k}}\widetilde{\Phi}_{H,k}\left(\bm{z}_{H,k}; \bm{z}_k^{(j)}\right)\right\|_2 \leq \frac{1}{2}\left\|\bm{W}_{\bm{u}}^H\right\|_{\infty} \left\|\nabla_{\bm{z}_k}\widetilde{\Phi}_k\left(\bm{z}_k; \bm{z}_k^{(j)}\right)\right\|_2$.
\end{itemize}
\end{proposition}
We leave the proof of this proposition to Appendix~\ref{app:proof_consistency}. 

{Furthermore, the parameters are chosen so that the coarse-grid solution of  Eqn.~\eqref{eqn:coarse_grid_surrogate_model}, obtained via the proximal point method, provides a descent direction for the fine-level problem. This is formalized in the following theorem.}
\begin{theorem}[Descent property]\label{prop:descent_direction}
    Given the weight vectors defined in~ Eqn.~\eqref{eqn:def_weights}, a proximal point update on the coarse surrogate \eqref{eqn:coarse_grid_surrogate_model}, initialized at $\bm{z}'_{H,k}= \bm I_h^H(\bm{W}_{\bm{z}}\odot\bm{z}'_k)$, provides a descent direction for the fine-grid surrogate \eqref{eqn:fine_grid_surrogate_model} at $\bm{z}'_k$. That is:
    \begin{equation*}       \nabla_{\bm{z}_k}\widetilde{\Phi}_k\left(\bm{z}_k'; \bm{z}_k^{(j)}\right)^*\left(\bm{I}_H^h\bm{e}_{H,k}\right)\leq0,
    \end{equation*}
    where $\bm{I}_H^h\bm{e}_{H,k} = \bm{I}_H^h\left(\widehat{\bm{z}_{H,k}}-\bm{z}_{H,k}'\right)$, and $\widehat{\bm{z}_{H,k}}$ solves:
    \begin{equation}\label{eqn:coarse_proximal}
            \widehat{\bm{z}_{H,k}} = \operatornamewithlimits{argmin}_{\bm{z}_{H,k}} \widetilde{\Phi}_{H,k}\left(\bm{z}_{H,k}; \bm{z}_k^{(j)}\right)+\frac{1}{2}\bm{u}_H^\top\left|\bm{z}_{H,k}-\bm{z}_{H,k}'\right|^2.
        \end{equation}
        Moreover,
        \begin{equation}\label{eqn:coarse_grid_rPIE_update}
            \bm{e}_{H,k} = -|\bm{Q}_H|^2\oslash\left(\bm{I}_h^H|\bm{Q}|^2\right)\oslash\left(\bm{u}_H+|\bm{Q}_H|^2\right)\odot\left(\bm{I}_h^H\nabla_{\bm{z}_k}\widetilde{\Phi}_k\left(\bm{z}_k'; \bm{z}_k^{(j)}\right)\right).
        \end{equation}
\end{theorem}
We leave the proof of this theorem to Appendix~\ref{app:proof_descent_direction}.

Whenever a model is represented on both fine and coarse grids, parameters fixed on the fine grid must have consistent counterparts on the coarse grid.  We show below that, with $\bm{W}_{\bm{z}},\bm{W}_{\mathcal{R}},\bm{W}_{\bm{u}}^H$ designed as in Eqn.~\eqref{eqn:def_weights}, the regularization vector can be transferred automatically between levels, requiring no manual tuning for the coarse-grid regularization vector $\bm{u}_H$. 

\begin{proposition}[Automatic regularization selection]\label{prop:automatic_regularization_selection}
Given a fine-grid regularization vector $\bm{u}$ and the weight vector $\bm{W}_{\bm{u}}^H$ defined in~ Eqn.~\eqref{eqn:def_weights}, by choosing
\begin{equation*}
    \bm{u}_H = \bm{W}_{\bm{u}}^H\odot \bm{I}_h^H\bm{u},
\end{equation*}
the induced coarse-grid regularization vector $\bm{u}_H$ ensures that, for each scanning region $k$, there exist mappings (with a slight abuse of notation)
\begin{equation*}
    \bm{e}_k:\mathbb{C}^{m^2}\to\mathbb{C}^{m^2}
    \quad\text{and}\quad
    \bm{e}_{H,k}:\mathbb{C}^{m^2}\to\mathbb{C}^{(m/2)^2},
\end{equation*}
such that:
\begin{itemize}
    \item The mappings $\bm{e}_k(\bm{g})$ and $\bm{e}_{H,k}(\bm{g})$ yield the fine-grid update at $\bm{z}_k'$ and the coarse-grid update at $\bm{z}_{H,k}' = \bm{I}_h^H(\bm{W}_{\bm{z}}\odot\bm{z}_k')$, respectively, for the fine-grid gradient $\bm{g} = \nabla_{\bm{z}_k}\widetilde{\Phi}_k(\bm{z}_k';\bm{z}_k^{(j)})$.
    
    \item For any index $r=1,\dots,(m/2)^2$, the coarse-grid update satisfies
    \begin{equation*}
        \max_{\|\bm{g}_1\|_\infty\le1}\left|\bm{e}_{H,k}(\bm{g}_1)\right|_r \leq \max_{\|\bm{g}_2\|_\infty\le1}\left|\bm{I}_h^H\bm{e}_k(\bm{g}_2)\right|_r. 
    \end{equation*}
    Thus, provided the fine‐grid update is stable under a properly tuned regularization vector, the maximum magnitude of the coarse‐grid update is automatically controlled.
\end{itemize}
\end{proposition}
We leave the proof of this proposition to Appendix~\ref{app:proof_automatic_regularization_selection}. 

\section{Numerical Examples}\label{sec:numerical_examples}
In this section, we evaluate MAGPIE and compare it with rPIE~\cite{rPIE} and L-BFGS~\cite{LBFGS}. For L-BFGS, we use history size $5$. To ensure a fair comparison with rPIE, we adopt its regularization as the fine-grid regularization in MAGPIE:
\begin{equation*}
  \bm{u}
  = \bm{u}^{\texttt{rPIE}}
  = \alpha\left(\|\bm{Q}\|_{\infty}^2\bm{J} - |\bm{Q}|^2\right).
\end{equation*}

The dominant per-iteration cost is evaluating the revised exit waves $\mathcal{R}_k\left(\bm{z}_k^{(j)}\right)$, which requires Fourier transforms for each scanning region. This term is computed once per scan position per iteration, and MAGPIE adds only element-wise operations and restriction/prolongation, whose cost is negligible by comparison. As the level deepens, the computational domain shrinks; at the deepest level, $\log_2(m)$, each scan update reduces to a scalar. Thus, MAGPIE has essentially the same per-iteration cost as rPIE and L-BFGS, while multilevel corrections improve convergence at little additional cost. This speedup is consistent with coarse-grid updates damping low-frequency error more efficiently by aggregating information over larger neighborhoods.

We first describe the experimental setup. Sections~\ref{sec:levels}, \ref{sec:noise}, and~\ref{sec:overlap_ratio} then assess performance with respect to multigrid depth, noise robustness, and overlap ratio. Section~\ref{sec:realistic_object} studies a simulated integrated-circuit object.

\begin{figure}[htbp]
    \centering
    \includegraphics[width=0.7\textwidth]{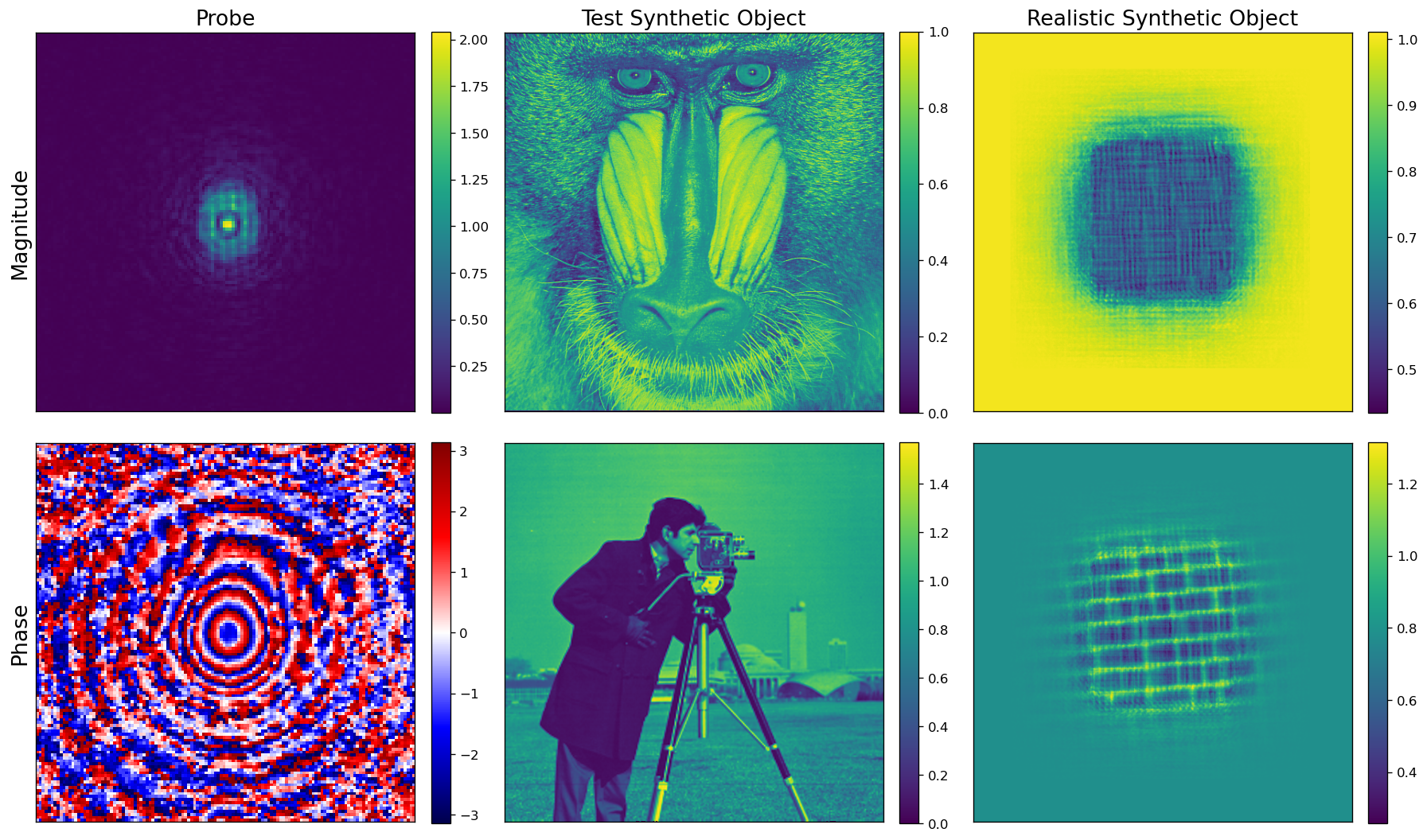}
    \caption{Magnitude (top) and phase (bottom) of the three inputs used in the numerical experiments: (1) the simulated Fresnel zone-plate probe; (2) the synthetic test object (Baboon magnitude and Cameraman phase); and (3) the realistic synthetic object inspired by integrated-circuit imaging.}
    \label{fig:all_mag_phase}
\end{figure}

\paragraph{Probe}
We use a simulated Fresnel zone plate as the probe (see~\cite{Ptychography}). The default probe size is $m=128$. Figure~\ref{fig:all_mag_phase} shows its magnitude and phase. The magnitude exhibits the characteristic circular structure of a zone plate, while the phase shows the corresponding alternating concentric pattern. Much of the probe outside the central region has very low magnitude, which may lead to ill-conditioning.

\paragraph{Object}
We consider two objects. The first is a synthetic test object formed from the Baboon image as magnitude and the Cameraman image as phase. The Baboon image is normalized to $[0,1]$, and the Cameraman image is scaled to $[0,\tfrac{\pi}{2}]$. The second is a more realistic synthetic object inspired by integrated-circuit imaging~\cite{integrated_circuit}. Both objects have resolution $n=512$ and are shown in Figure~\ref{fig:all_mag_phase}.

\paragraph{Overlap ratio}
In our experiments, the overlap ratio is set to 50\% by default, i.e., $\texttt{overlap\_ratio}=0.5$.

\paragraph{Noise level}
To mimic photon-counting statistics, we add Poisson noise to each simulated intensity. Given the noiseless value $\widetilde{\bm{d}}_k$, the noisy measurement is
\begin{equation*}
    \bm{d}_k = \eta\,\mathrm{Poisson}\!\left(\widetilde{\bm{d}}_k/\eta\right),
\end{equation*}
where the dimensionless parameter $\eta$ controls the noise level: larger $\eta$ lowers the mean photon count and hence increases the relative noise. Unless otherwise stated, we use $\eta=0.05$.

\paragraph{Stopping criterion}
For all algorithms, we stop when
\begin{equation*}
    \frac{1}{Nm}\sum_{k=1}^{N}\left\|\nabla_{\bm z_k}\Phi_k(\bm z_k)\right\|_2 < \texttt{tol}.
\end{equation*}

\paragraph{Quality metrics}
We evaluate reconstructions using two metrics: residual and error. For a reconstruction $\bm{z}$, the residual is the loss value $\Phi(\bm{z})$. The error is the $\ell^2$ difference between the reconstructed and true magnitudes, i.e.,
\begin{equation*}
    \||\bm{z}| - |\bm{z}^*|\|_2.
\end{equation*}

{In all convergence plots below, each method is run for at most 3000 iterations. A hollow circle marks the first iteration satisfying the stopping criterion, and dashed segments indicate subsequent iterations. Unless otherwise stated, reconstruction images are taken from the first iterate satisfying the stopping criterion.}
\subsection{Effect of multigrid levels on performance}\label{sec:levels}
We study how the number of coarse levels affects MAGPIE. For $l=1,\ldots,\log_2(m)$, we denote the corresponding method by MAGPIE\_$l$. Since $m=128$, the deepest hierarchy has $\log_2(m)=7$ levels. We use $\texttt{tol}=10^{-4}$ and $\alpha=0.01$ throughout.

\begin{figure}[htbp]
    \centering
    \includegraphics[width=0.8\textwidth]{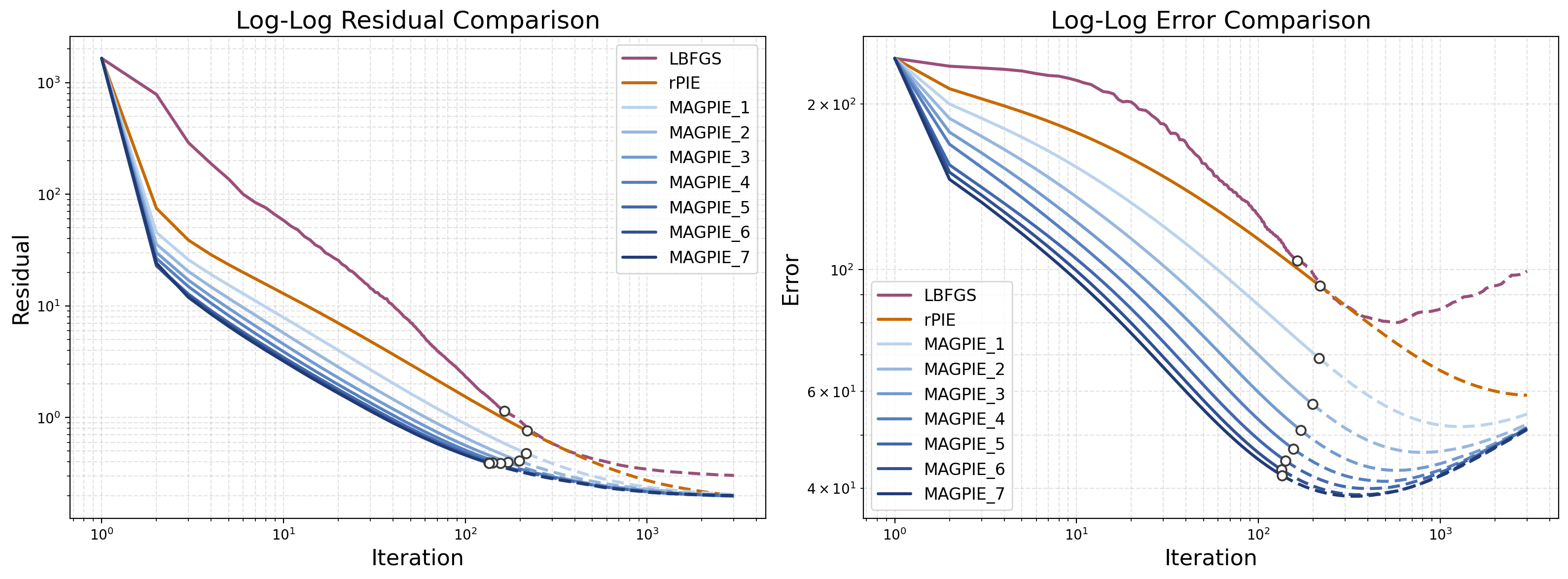}
    \caption{Log-log plots of residuals (top) and errors (bottom) for L-BFGS, rPIE, and MAGPIE\_$l$ with $l=1,\ldots,\log_2(m)$, applied to a synthetic object ($n=512$) with a probe ($m=128$), noise level $\eta=0.05$, $\texttt{overlap\_ratio}=0.5$, $\alpha=0.01$, and $\texttt{tol}=10^{-4}$.}
    \label{fig:metric_1}
\end{figure}

{Figure~\ref{fig:metric_1} shows residuals and errors on a log-log scale. Deeper multilevel variants consistently satisfy the stopping criterion earlier and achieve smaller reconstruction errors. Over the full 3000-iteration budget, the deeper MAGPIE variants remain more accurate than both rPIE and L-BFGS while also attaining the smallest residuals. Accordingly, in the remaining experiments we use the deepest version, $l=\log_2(m)$.}

\begin{figure}[htbp]
    \centering
    \includegraphics[width=0.8\textwidth]{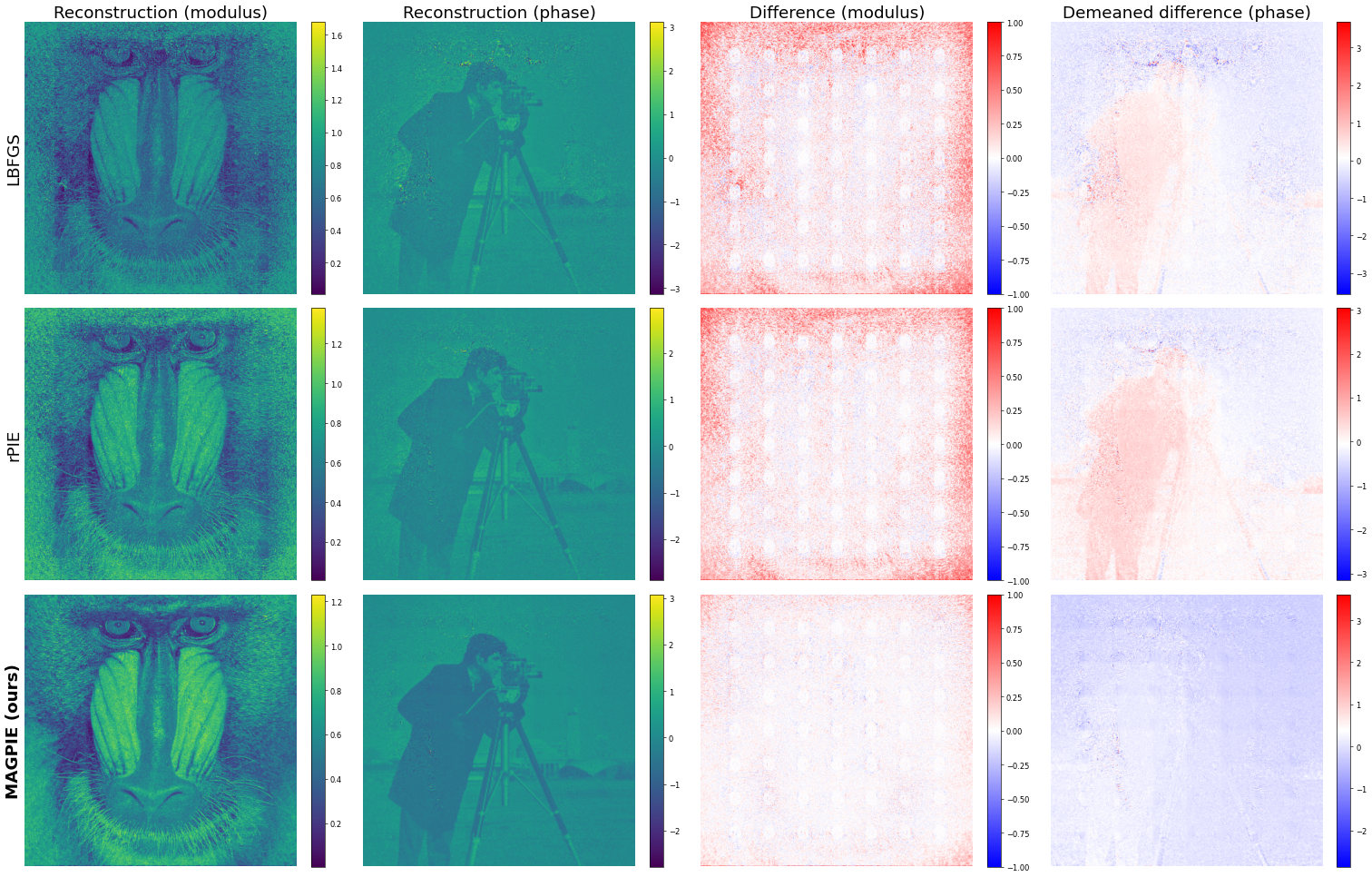}
    \caption{Reconstructions of L-BFGS, rPIE, and MAGPIE at level $\log_2(m)$ for a synthetic object ($n=512$) with a probe ($m=128$), noise level $\eta=0.05$, $\texttt{overlap\_ratio}=0.5$, $\alpha=0.01$, and $\texttt{tol}=10^{-4}$.}
    \label{fig:reconstruction_1}
\end{figure}

Figure~\ref{fig:reconstruction_1} shows the corresponding reconstructions. MAGPIE yields the best reconstruction quality, with smaller and smoother magnitude and phase errors than L-BFGS and rPIE.

\subsection{Stability under noise}\label{sec:noise}
We next assess robustness to noise by comparing MAGPIE with L-BFGS and rPIE at four Poisson noise levels: $\eta=0.05, 0.1, 0.2,$ and $0.4$. Throughout, we set $\texttt{tol}=10^{-4}$ and $\alpha=0.025$.

\begin{figure}[htbp]
    \centering
    \includegraphics[width=0.8\textwidth]{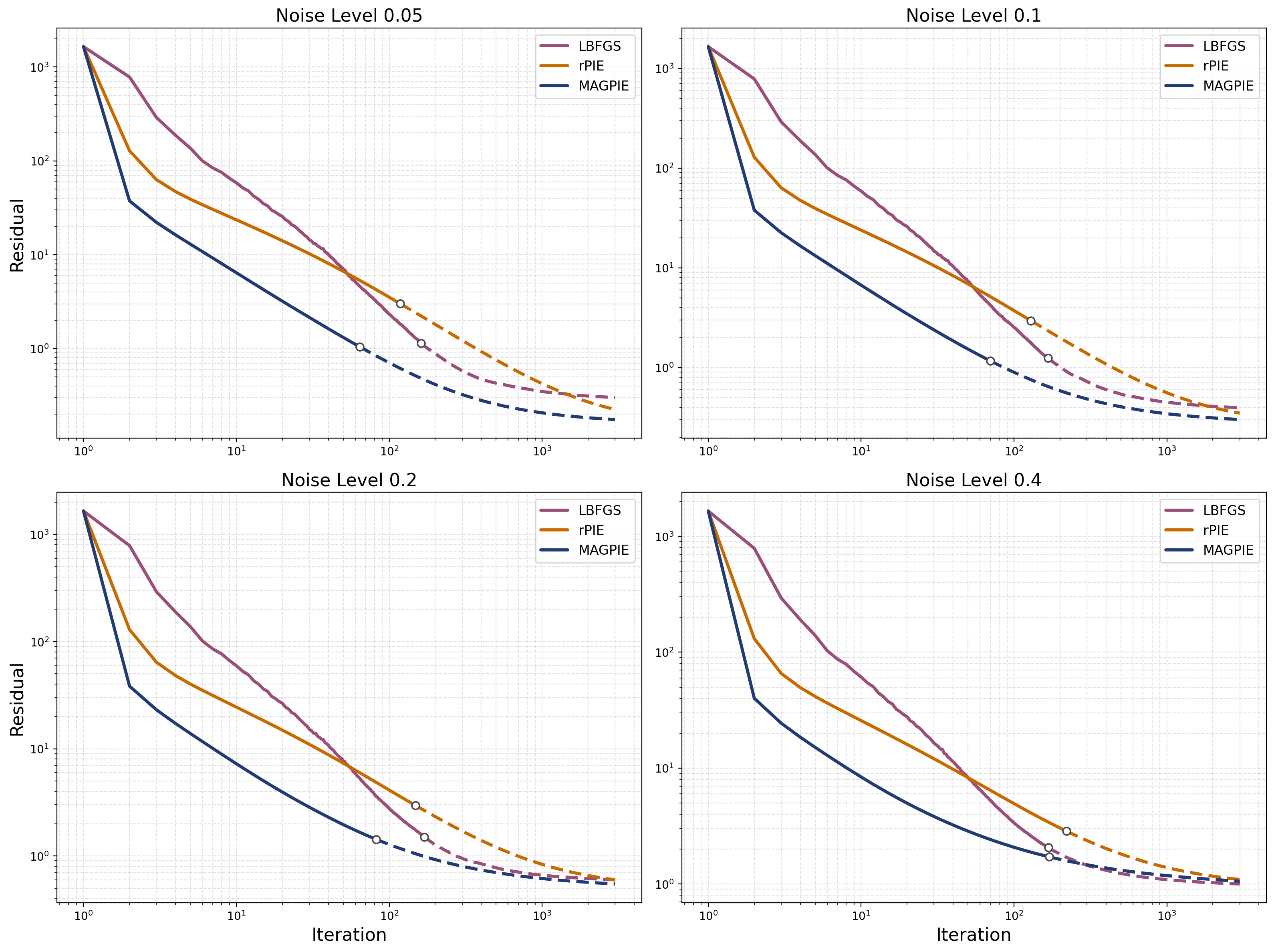}
    \caption{Log-log plots of residuals for L-BFGS, rPIE, and MAGPIE applied to a synthetic object ($n=512$) with a probe ($m=128$), noise level $\eta=0.05, 0.1, 0.2,$ and $0.4$, $\texttt{overlap\_ratio}=0.5$, $\alpha=0.025$, and $\texttt{tol}=10^{-4}$.}
    \label{fig:residual_2}
\end{figure}

\begin{figure}[htbp]
    \centering
    \includegraphics[width=0.8\textwidth]{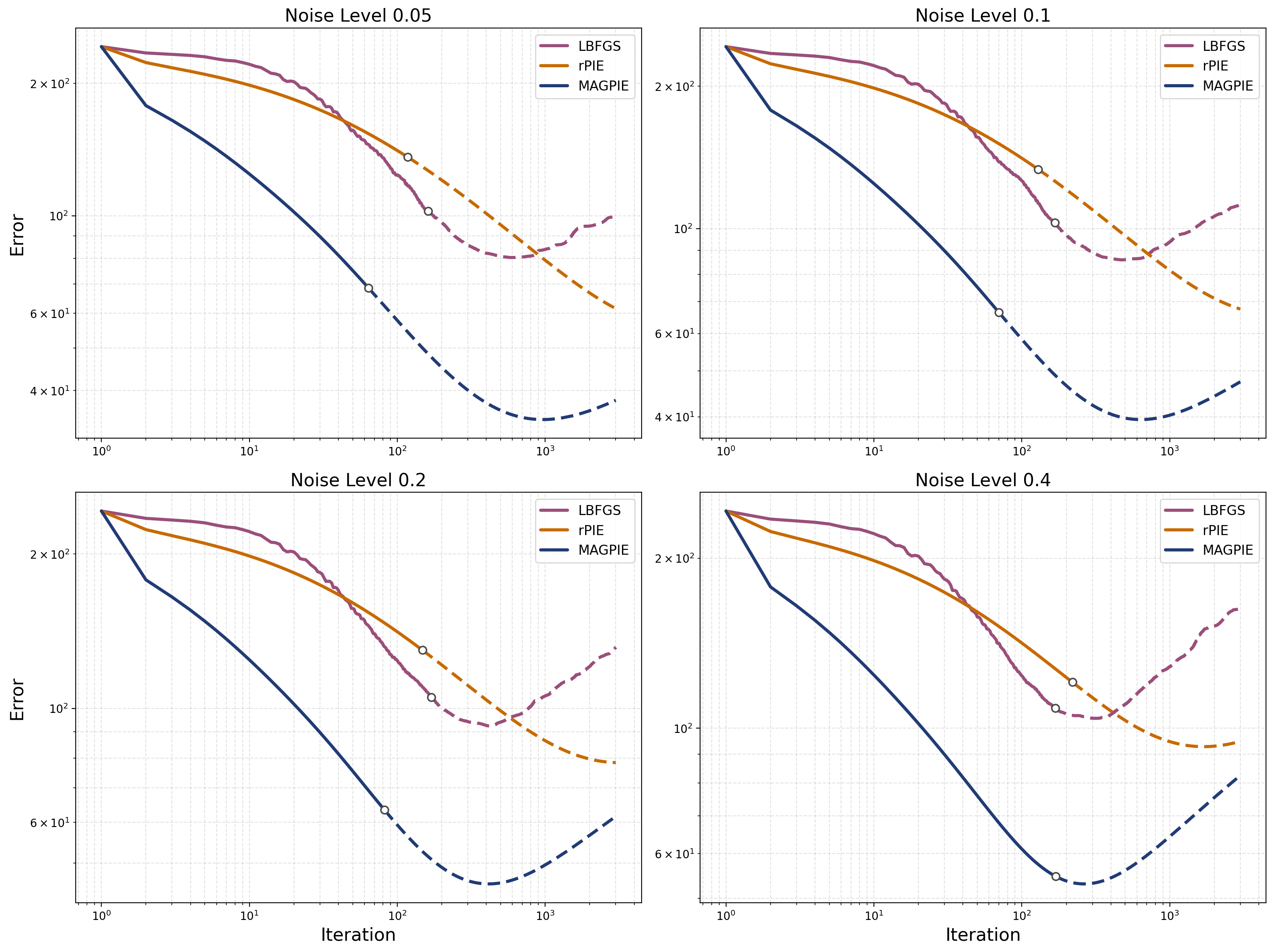}
    \caption{Log-log plots of errors for L-BFGS, rPIE, and MAGPIE applied to a synthetic object ($n=512$) with a probe ($m=128$), noise level $\eta=0.05, 0.1, 0.2,$ and $0.4$, $\texttt{overlap\_ratio}=0.5$, $\alpha=0.025$, and $\texttt{tol}=10^{-4}$.}
    \label{fig:error_2}
\end{figure}

{Figures~\ref{fig:residual_2} and~\ref{fig:error_2} show the residuals and errors on a log-log scale. Across all noise levels, MAGPIE reaches the stopping criterion earlier and attains the smallest reconstruction error in the practically relevant regime. Over the full 3000-iteration budget, it also remains the most accurate. The slight upturns in some dashed late-iteration error curves are consistent with fitting measurement noise after the stopping point rather than improved reconstruction quality.}

\begin{figure}[htbp]
    \centering
    \includegraphics[width=0.8\textwidth]{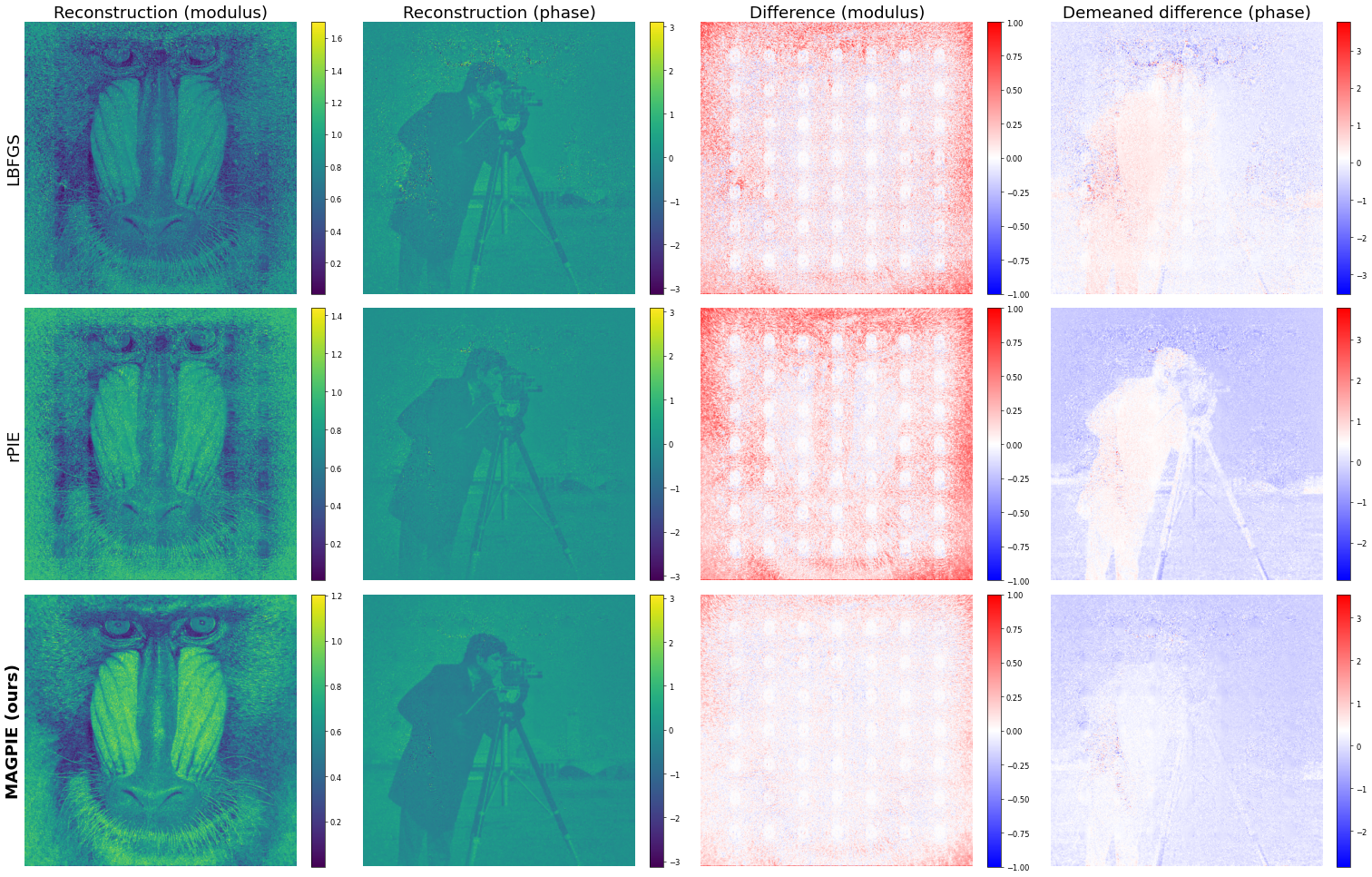}
    \caption{Reconstructions of L-BFGS, rPIE, and MAGPIE applied to a synthetic object ($n=512$) with a probe ($m=128$), noise level $\eta=0.4$, $\texttt{overlap\_ratio}=0.5$, $\alpha=0.025$, and $\texttt{tol}=10^{-4}$.}
    \label{fig:reconstruction_2}
\end{figure}

Figure~\ref{fig:reconstruction_2} shows the reconstructions at noise level $\eta=0.4$. MAGPIE again gives the best result, with smaller and smoother errors than L-BFGS and rPIE.

\subsection{Overlap ratio}\label{sec:overlap_ratio}
We next study the effect of probe overlap by comparing MAGPIE with rPIE and L-BFGS at two nondefault overlap ratios: $\texttt{overlap\_ratio}=0.25$, with $\alpha=0.008$ and $\texttt{tol}=10^{-5}$, and $\texttt{overlap\_ratio}=0.75$, with $\alpha=0.04$ and $\texttt{tol}=10^{-4}$.

\begin{figure}[htbp]
    \centering
    \includegraphics[width=0.8\textwidth]{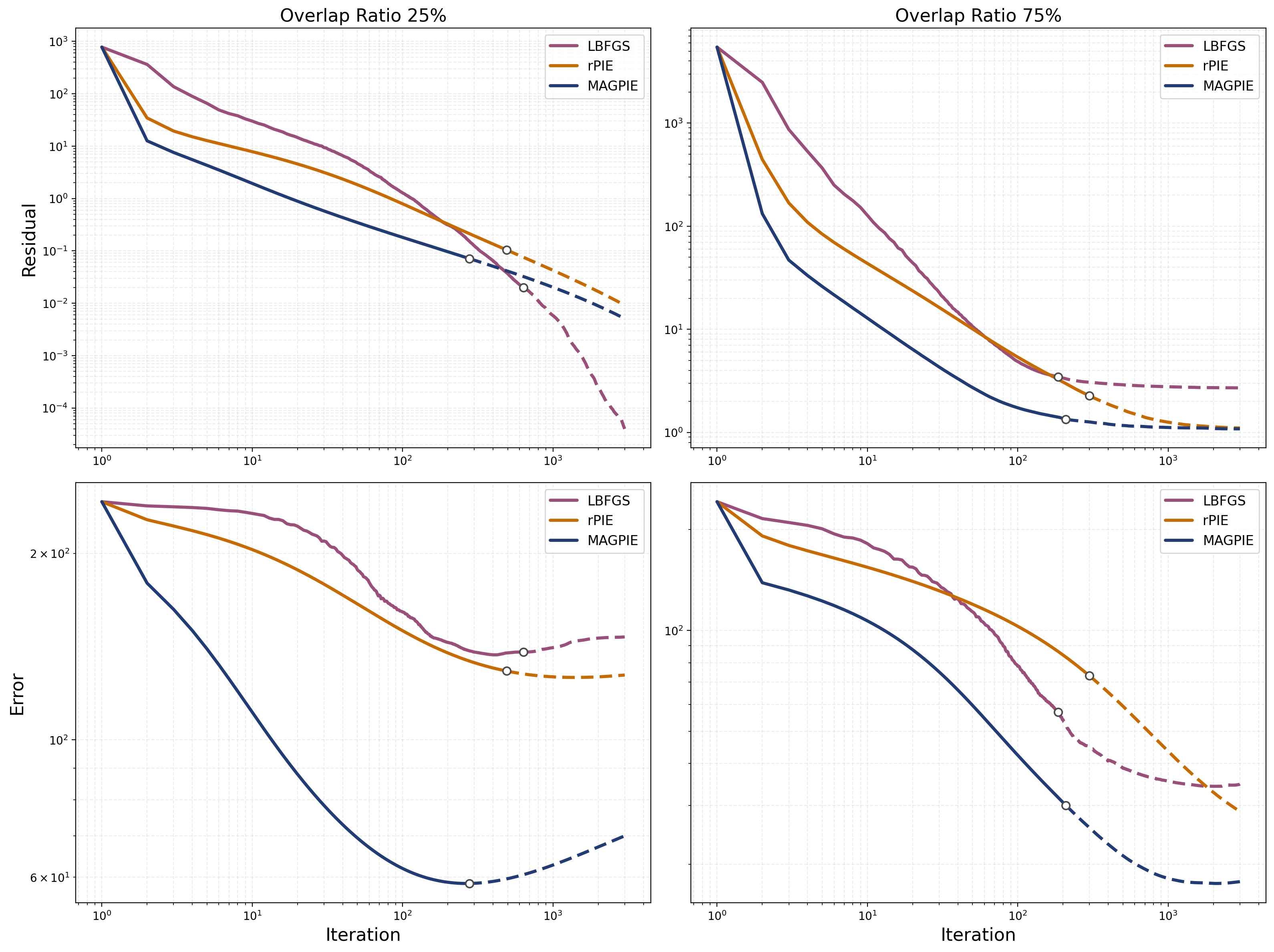}
    \caption{Log-log plots of residuals and errors for L-BFGS, rPIE, and MAGPIE applied to a synthetic object ($n=512$) with probe size $m=128$ and noise level $\eta=0.05$. The first column corresponds to $\texttt{overlap\_ratio}=0.25$, $\alpha=0.008$, and $\texttt{tol}=10^{-5}$; the second corresponds to $\texttt{overlap\_ratio}=0.75$, $\alpha=0.04$, and $\texttt{tol}=10^{-4}$.}
    \label{fig:metric_3}
\end{figure}

\begin{figure}[htbp]
    \centering
    \includegraphics[width=0.8\textwidth]{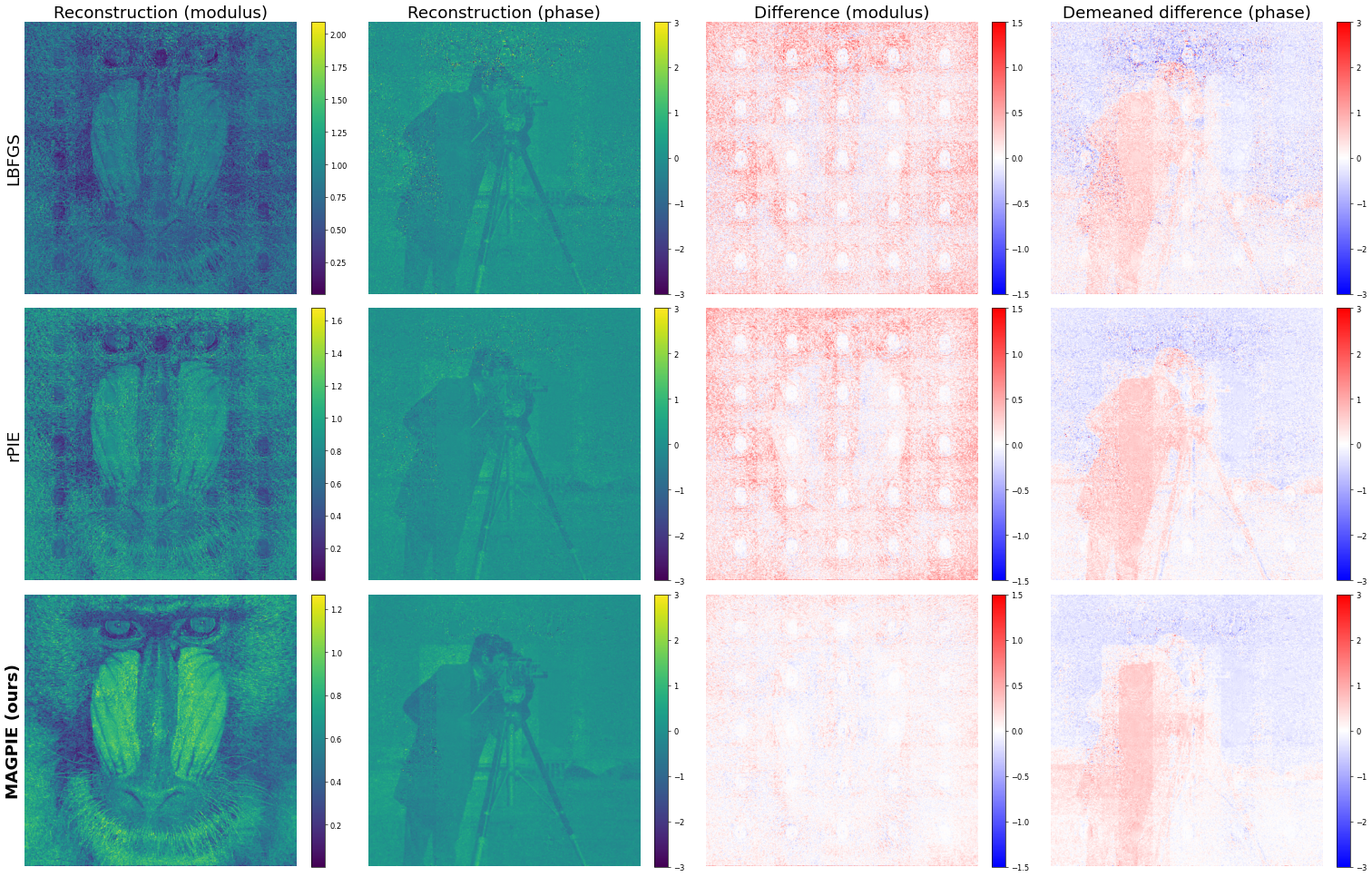}
    \caption{Reconstructions of L-BFGS, rPIE, and MAGPIE applied to a synthetic object ($n=512$) with a probe ($m=128$), noise level $\eta=0.05$, $\texttt{overlap\_ratio}=0.25$, $\alpha=0.008$, and $\texttt{tol}=10^{-5}$.}
    \label{fig:reconstruction_3_25}
\end{figure}

\begin{figure}[htbp]
    \centering
    \includegraphics[width=0.8\textwidth]{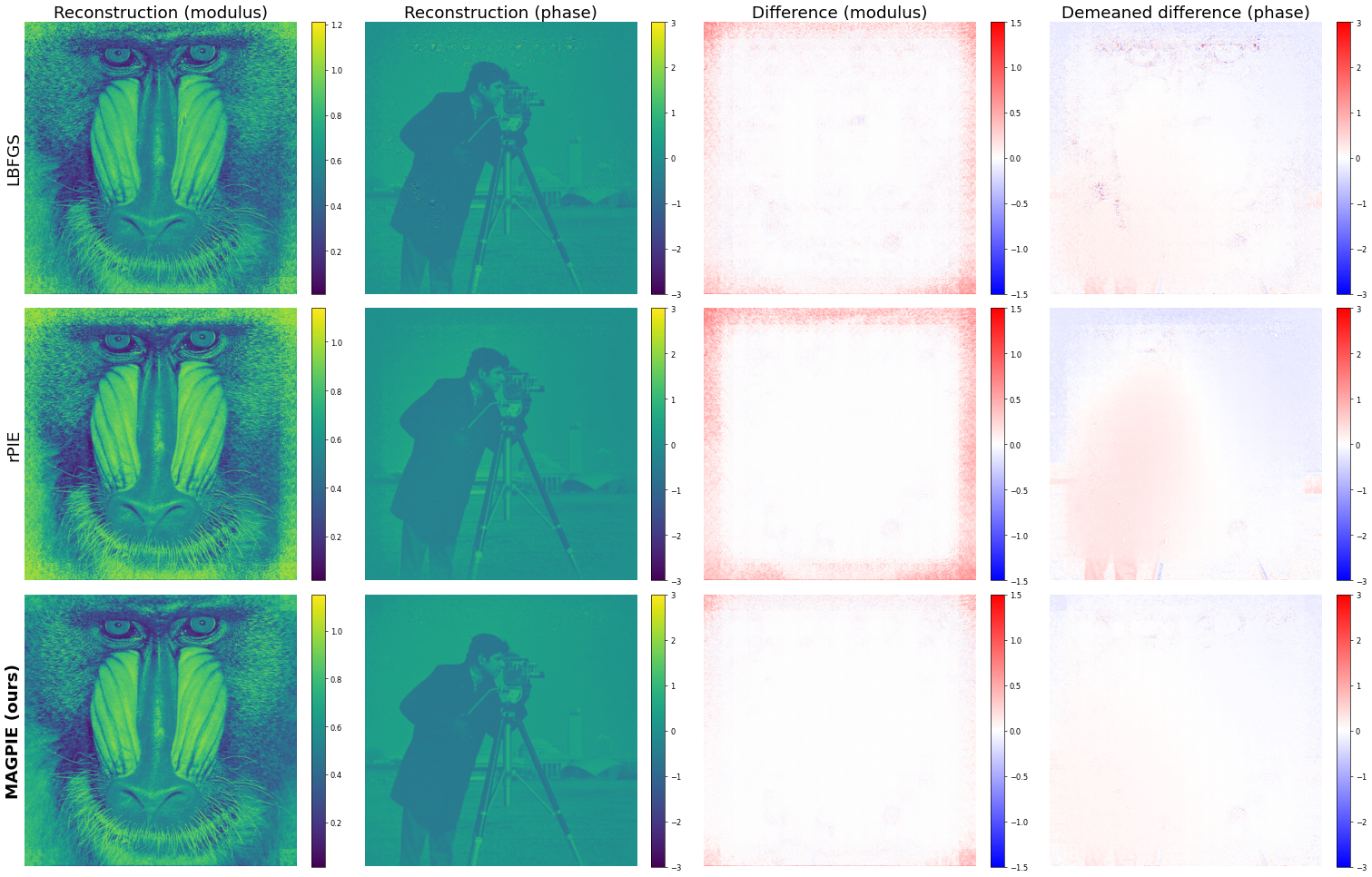}
    \caption{Reconstructions of L-BFGS, rPIE, and MAGPIE applied to a synthetic object ($n=512$) with a probe ($m=128$), noise level $\eta=0.05$, $\texttt{overlap\_ratio}=0.75$, $\alpha=0.04$, and $\texttt{tol}=10^{-4}$.}
    \label{fig:reconstruction_3_75}
\end{figure}

Figure~\ref{fig:metric_3} shows residuals and errors for both cases. At $\texttt{overlap\_ratio}=0.25$, L-BFGS eventually attains the smallest residual after the stopping point, but its reconstruction error remains substantially larger than that of MAGPIE, indicating that the lower residual does not yield better reconstruction quality. This is consistent with the reduced data redundancy and weaker constraints at low overlap. In this regime, MAGPIE reaches low error much earlier and remains the most accurate method over the extended horizon. Reconstructions in Figures~\ref{fig:reconstruction_3_25} and~\ref{fig:reconstruction_3_75} are shown at the stopping point.

\subsection{Realistic synthetic object}\label{sec:realistic_object}
We repeat the comparison on a more realistic synthetic object whose entries are mostly zero near the boundary. This test shows that MAGPIE’s performance gain is not merely due to improved reconstruction near the boundary, where fewer scan regions overlap.

We consider two cases: $\texttt{overlap\_ratio}=0.50$, with $\alpha=0.03$ and $\texttt{tol}=3\times10^{-5}$, and $\texttt{overlap\_ratio}=0.75$, with $\alpha=0.06$ and $\texttt{tol}=6\times10^{-5}$.

\begin{figure}[htbp]
    \centering
    \includegraphics[width=0.8\textwidth]{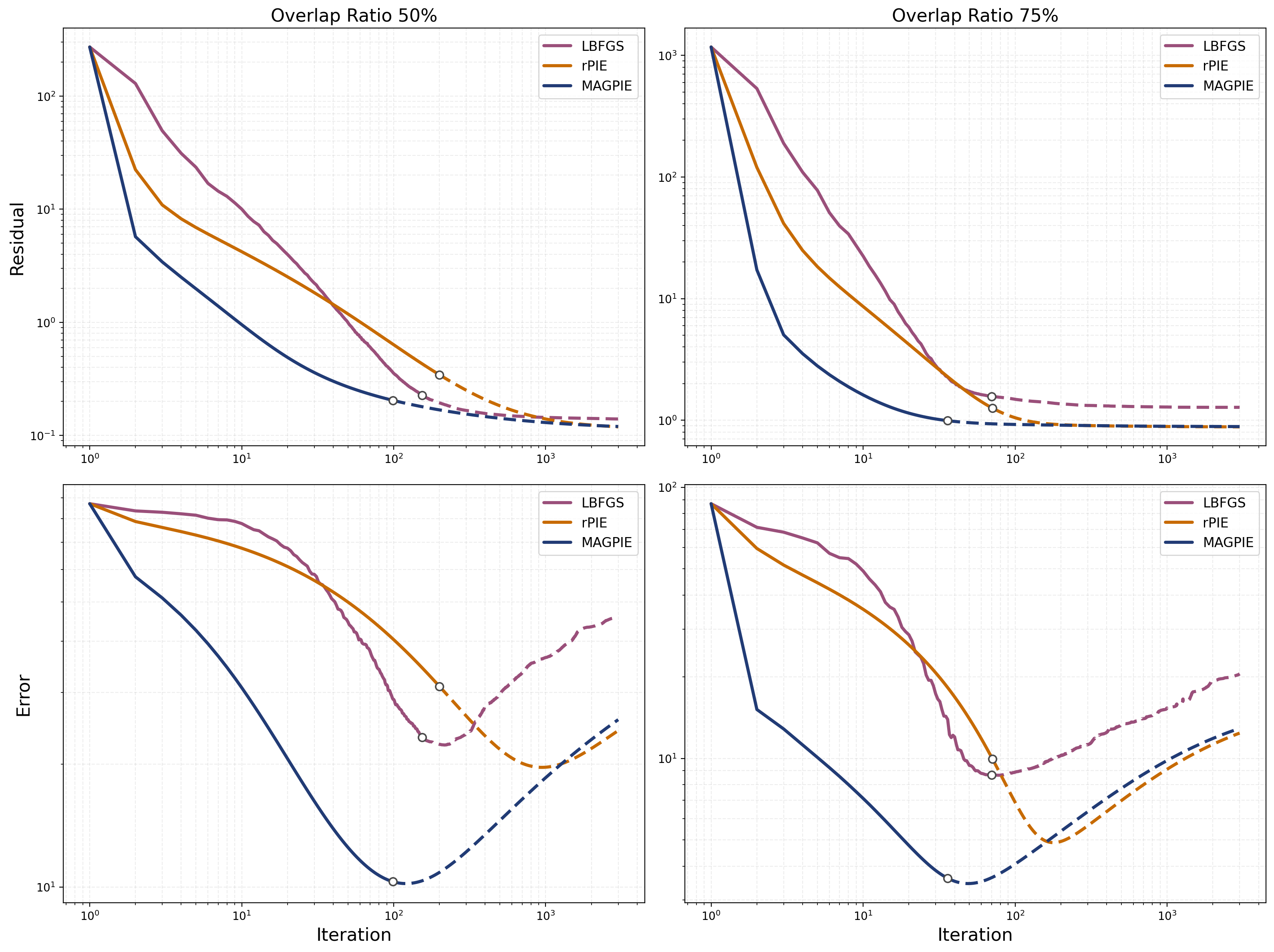}
    \caption{Log-log plots of residuals and errors for L-BFGS, rPIE, and MAGPIE applied to a realistic synthetic object ($n=512$) with probe size $m=128$ and noise level $\eta=0.05$. The first column corresponds to $\texttt{overlap\_ratio}=0.50$, $\alpha=0.03$, and $\texttt{tol}=3\times10^{-5}$; the second corresponds to $\texttt{overlap\_ratio}=0.75$, $\alpha=0.06$, and $\texttt{tol}=6\times10^{-5}$.}
    \label{fig:metric_4}
\end{figure}

\begin{figure}[htbp]
    \centering
    \includegraphics[width=0.8\textwidth]{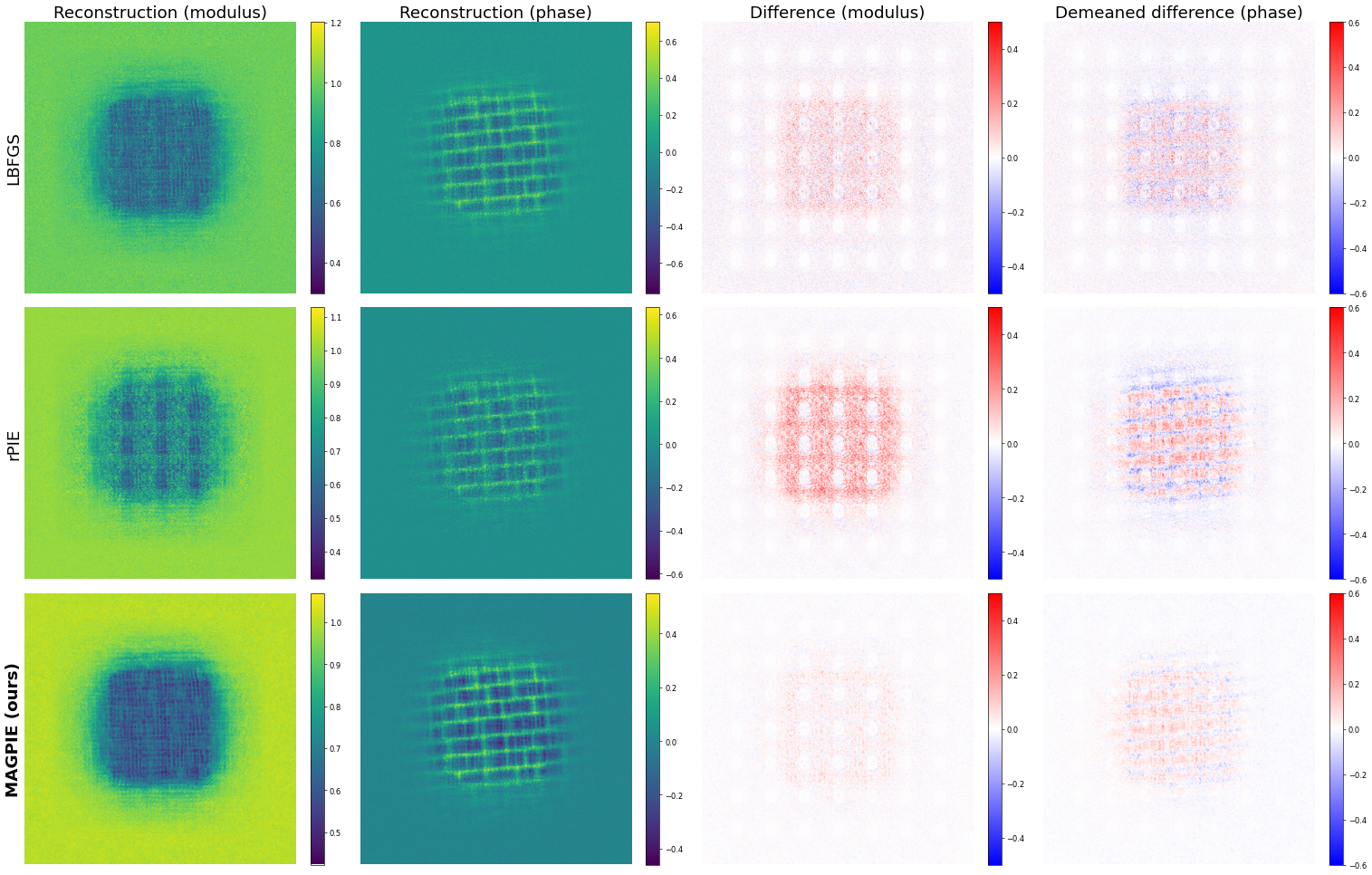}
    \caption{Reconstructions of L-BFGS, rPIE, and MAGPIE applied to a realistic synthetic object ($n=512$) with a probe ($m=128$), noise level $\eta=0.05$, $\texttt{overlap\_ratio}=0.50$, $\alpha=0.03$, and $\texttt{tol}=3\times10^{-5}$.}
    \label{fig:reconstruction_4_50}
\end{figure}

\begin{figure}[htbp]
    \centering
    \includegraphics[width=0.8\textwidth]{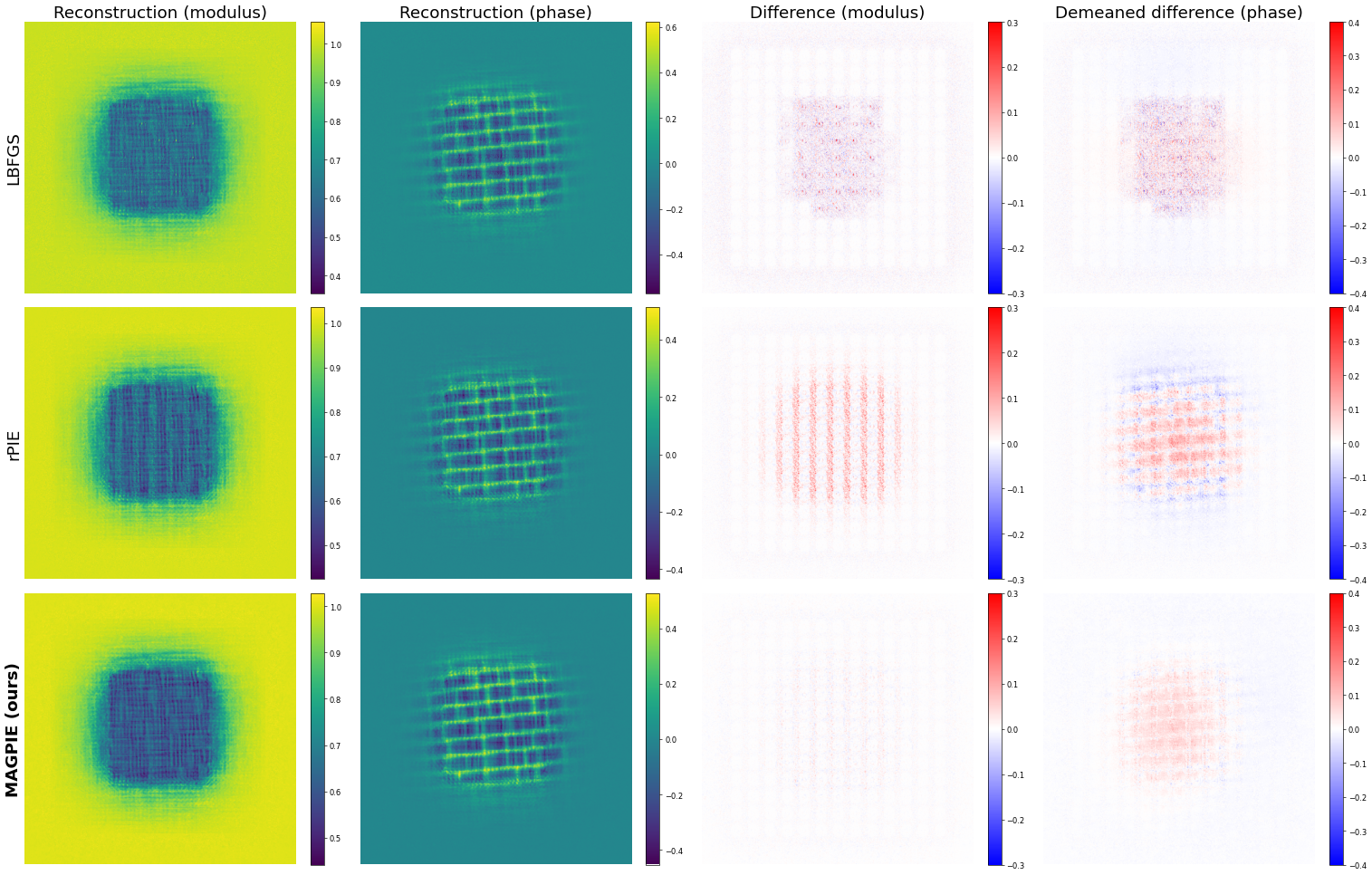}
    \caption{Reconstructions of L-BFGS, rPIE, and MAGPIE applied to a realistic synthetic object ($n=512$) with a probe ($m=128$), noise level $\eta=0.05$, $\texttt{overlap\_ratio}=0.75$, $\alpha=0.06$, and $\texttt{tol}=6\times10^{-5}$.}
    \label{fig:reconstruction_4_75}
\end{figure}

{Figure~\ref{fig:metric_4} shows the residuals and errors for both cases. In both overlap settings, MAGPIE reaches a low-error reconstruction substantially earlier than rPIE and L-BFGS. With continued iteration, rPIE can partially catch up in error after the stopping point, whereas L-BFGS plateaus or deteriorates. This is consistent with late-stage fitting of noisy measurements after MAGPIE has already reached a high-quality iterate. Reconstructions in Figures~\ref{fig:reconstruction_4_50} and~\ref{fig:reconstruction_4_75} are shown at the stopping point.}

\section{Conclusion}
We have presented MAGPIE, a stochastic multigrid minimization framework for ptychographic phase retrieval that accelerates PIE-type updates via rigorously constructed coarse-grid surrogate models. The surrogate guarantees monotonic descent of the original nonconvex objective with sublinear convergence, embeds a multilevel hierarchy in which coarse-grid regularization vectors are automatically selected by fine–coarse consistency, and retains the standard rPIE parameter set while significantly improving both convergence speed and reconstruction fidelity.

Looking ahead, we will extend the theoretical framework to blind ptychography and integrate a minibatch strategy for GPU-accelerated real‐time feedback. Moreover, the multigrid‐surrogate paradigm developed here is readily adaptable to other large‐scale phase retrieval and coherent diffractive imaging applications.

\section{Acknowledgments}
Z.D. was supported by the U.S. Department of Energy, Office of Science, Advanced Scientific Computing Research, under Contract DE-AC02-06CH11357. Q.L. acknowledges support from NSF-DMS-2308440, and B.Z. acknowledges support from NSF-DMS-2012292.
\bibliographystyle{unsrt}
\bibliography{references} 

@article{SMM_graph,
author = {Mairal, Julien},
title = {Incremental Majorization-Minimization Optimization with Application to Large-Scale Machine Learning},
journal = {SIAM J. Optim.},
volume = {25},
number = {2},
pages = {829--855},
year = {2015},
doi = {10.1137/140957639},
URL = {https://doi.org/10.1137/140957639},
eprint = {https://doi.org/10.1137/140957639},
}

@article{mginverse1,
doi = {10.1088/0266-5611/17/4/313},
url = {https://doi.org/10.1088/0266-5611/17/4/313},
year = {2001},
month = {aug},
publisher = {},
volume = {17},
number = {4},
pages = {767--788},
author = {Barbara Kaltenbacher},
title = {On the regularizing properties of a full multigrid method for ill-posed problems},
journal = {Inverse Problems}}

@INPROCEEDINGS{MM_phaseretrieval,
  author={Qiu, Tianyu and Babu, Prabhu and Palomar, Daniel P.},
  booktitle={2015 49th Asilomar Conference on Signals, Systems and Computers}, 
  title={{PRIME}: Phase retrieval via majorization-minimization technique}, 
  year={2015},
  volume={},
  number={},
  pages={1681--1685},
  doi={10.1109/ACSSC.2015.7421435}}

@book{multigrid_book,
  title = {A Multigrid Tutorial,  Second Edition},
  ISBN = {9780898719505},
  url = {http://dx.doi.org/10.1137/1.9780898719505},
  DOI = {10.1137/1.9780898719505},
  publisher = {Society for Industrial and Applied Mathematics},
  author = {Briggs,  William L. and Henson,  Van Emden and McCormick,  Steve F.},
  year = {2000},
  month = Jan 
}

@article{mgopt_1, 
author  = {Di, Zichao and Emelianenko, Maria and Nash, Stephen},
title   = {Truncated {Newton}-Based Multigrid Algorithm for Centroidal {Voronoi} Diagram Calculation},
volume = {5},
ISSN = {1004-8979},
url = {http://dx.doi.org/10.4208/nmtma.2012.m1046},
DOI = {10.4208/nmtma.2012.m1046},
number = {2},
journal = {Numer. Math. Theory Methods Appl.},
publisher = {Global Science Press},
year = {2012},
month = Jan,
pages = {242--259}}

@incollection{Rodenburg2019,
  title = {Ptychography},
  ISBN = {9783030000691},
  ISSN = {2522-8706},
  url = {http://dx.doi.org/10.1007/978-3-030-00069-1_17},
  DOI = {10.1007/978-3-030-00069-1_17},
  booktitle = {Springer Handbook of Microscopy},
  publisher = {Springer International Publishing},
  author = {Rodenburg,  John and Maiden,  Andrew},
  year = {2019},
  pages = {819--904}
}

@article{crystallography,
  title = {Ptychographic Imaging of Branched Colloidal Nanocrystals Embedded in Free-Standing Thick Polystyrene Films},
  volume = {6},
  ISSN = {2045-2322},
  url = {http://dx.doi.org/10.1038/srep19397},
  DOI = {10.1038/srep19397},
  number = {1},
  journal = {Scientific Reports},
  publisher = {Springer Science and Business Media LLC},
  author = {De Caro,  Liberato and Altamura,  Davide and Arciniegas,  Milena and Siliqi,  Dritan and Kim,  Mee R. and Sibillano,  Teresa and Manna,  Liberato and Giannini,  Cinzia},
  year = {2016},
  month = Jan 
}

@article{biology_2,
  title = {Dark-field {X}-ray ptychography: Towards high-resolution imaging of thick and unstained biological specimens},
  volume = {6},
  ISSN = {2045-2322},
  url = {http://dx.doi.org/10.1038/srep35060},
  DOI = {10.1038/srep35060},
  number = {1},
  journal = {Scientific Reports},
  publisher = {Springer Science and Business Media LLC},
  author = {Suzuki,  Akihiro and Shimomura,  Kei and Hirose,  Makoto and Burdet,  Nicolas and Takahashi,  Yukio},
  year = {2016},
  month = Oct 
}

@article{biology_1,
  title = {Ptychography – a label free,  high-contrast imaging technique for live cells using quantitative phase information},
  volume = {3},
  ISSN = {2045-2322},
  url = {http://dx.doi.org/10.1038/srep02369},
  DOI = {10.1038/srep02369},
  number = {1},
  journal = {Scientific Reports},
  publisher = {Springer Science and Business Media LLC},
  author = {Marrison,  Joanne and R\"{a}ty,  Lotta and Marriott,  Poppy and O’Toole,  Peter},
  year = {2013},
  month = Aug 
}

@article{integrated_circuit_1,
  title = {Achieving high spatial resolution in a large field-of-view using lensless {X}-ray imaging},
  volume = {119},
  ISSN = {1077-3118},
  url = {http://dx.doi.org/10.1063/5.0067197},
  DOI = {10.1063/5.0067197},
  number = {12},
  journal = {Applied Physics Letters},
  publisher = {AIP Publishing},
  author = {Jiang,  Yi and Deng,  Junjing and Yao,  Yudong and Klug,  Jeffrey A. and Mashrafi,  Sheikh and Roehrig,  Christian and Preissner,  Curt and Marin,  Fabricio S. and Cai,  Zhonghou and Lai,  Barry and Vogt,  Stefan},
  year = {2021},
  month = {September},
}

@article{materials_science_2,
  title = {High-resolution chemical imaging of gold nanoparticles using hard {X}-ray ptychography},
  volume = {102},
  ISSN = {1077-3118},
  url = {http://dx.doi.org/10.1063/1.4807020},
  DOI = {10.1063/1.4807020},
  number = {20},
  journal = {Applied Physics Letters},
  publisher = {AIP Publishing},
  author = {Hoppe,  R. and Reinhardt,  J. and Hofmann,  G. and Patommel,  J. and Grunwaldt,  J.-D. and Damsgaard,  C. D. and Wellenreuther,  G. and Falkenberg,  G. and Schroer,  C. G.},
  year = {2013},
  month = May 
}

@article{integrated_circuit,
  title = {High-resolution non-destructive three-dimensional imaging of integrated circuits},
  volume = {543},
  ISSN = {1476-4687},
  url = {http://dx.doi.org/10.1038/nature21698},
  DOI = {10.1038/nature21698},
  number = {7645},
  journal = {Nature},
  publisher = {Springer Science and Business Media LLC},
  author = {Holler,  Mirko and Guizar-Sicairos,  Manuel and Tsai,  Esther H. R. and Dinapoli,  Roberto and M\"{u}ller,  Elisabeth and Bunk,  Oliver and Raabe,  J\"{o}rg and Aeppli,  Gabriel},
  year = {2017},
  month = Mar,
  pages = {402--406}
}

@incollection{Ptychography,
  title = {Ptychography and Related Diffractive Imaging Methods},
  ISBN = {9780123742179},
  ISSN = {1076-5670},
  url = {http://dx.doi.org/10.1016/S1076-5670(07)00003-1},
  DOI = {10.1016/s1076-5670(07)00003-1},
  booktitle = {Advances in Imaging and Electron Physics},
  publisher = {Elsevier},
  author = {Rodenburg,  J.M.},
  year = {2008},
  pages = {87--184}
}

@article{LBFGS,
  title = {On the limited memory {BFGS} method for large scale optimization},
  volume = {45},
  ISSN = {1436-4646},
  url = {http://dx.doi.org/10.1007/BF01589116},
  DOI = {10.1007/bf01589116},
  number = {1-3},
  journal = {Math. Program.},
  publisher = {Springer Science and Business Media LLC},
  author = {Liu,  Dong C. and Nocedal,  Jorge},
  year = {1989},
  month = Aug,
  pages = {503--528}
}

@article{MGOPT,
  title = {Properties of a class of multilevel optimization algorithms for equality-constrained problems},
  volume = {29},
  ISSN = {1029-4937},
  url = {http://dx.doi.org/10.1080/10556788.2012.759571},
  DOI = {10.1080/10556788.2012.759571},
  number = {1},
  journal = {Optim. Methods Softw.},
  publisher = {Informa UK Limited},
  author = {Nash,  Stephen G.},
  year = {2014},
  month = Mar,
  pages = {137--159}
}

@article{mm_alg_Lange,
  title = {Optimization Transfer Using Surrogate Objective Functions},
  volume = {9},
  ISSN = {1537-2715},
  url = {http://dx.doi.org/10.1080/10618600.2000.10474858},
  DOI = {10.1080/10618600.2000.10474858},
  number = {1},
  journal = {J. Comput. Graph. Statist.},
  publisher = {Informa UK Limited},
  author = {Lange,  Kenneth and Hunter,  David R. and Yang,  Ilsoon},
  year = {2000},
  month = Mar,
  pages = {1--20}
}

@incollection{mm_convergence,
  title = {Nonconvex Optimization via {MM} Algorithms: Convergence Theory},
  ISBN = {9781118445112},
  url = {http://dx.doi.org/10.1002/9781118445112.stat08295},
  DOI = {10.1002/9781118445112.stat08295},
  booktitle = {Wiley StatsRef: Statistics Reference Online},
  publisher = {John Wiley \& Sons, Ltd},
  author = {Lange,  Kenneth and Won,  Joong‐Ho and Landeros,  Alfonso and Zhou,  Hua},
  year = {2021},
  month = Mar,
  pages = {1--22}
}

@article{MGOPT_ptycho,
  title = {Multigrid Optimization for Large-Scale Ptychographic Phase Retrieval},
  volume = {13},
  ISSN = {1936-4954},
  url = {http://dx.doi.org/10.1137/18M1223915},
  DOI = {10.1137/18m1223915},
  number = {1},
  journal = {SIAM J. Imaging Sci.},
  publisher = {Society for Industrial & Applied Mathematics (SIAM)},
  author = {Wu Fung, Samy and Di, Zichao Wendy},
  year = {2020},
  month = Jan,
  pages = {214--233}
}

@article{Nash01012000,
  title = {A multigrid approach to discretized optimization problems},
  volume = {14},
  ISSN = {1029-4937},
  url = {http://dx.doi.org/10.1080/10556780008805795},
  DOI = {10.1080/10556780008805795},
  number = {1-2},
  journal = {Optim. Methods Softw.},
  publisher = {Informa UK Limited},
  author = {Nash,  Stephen G.},
  year = {2000},
  month = Jan,
  pages = {99--116}
}

@article{PIE,
  title = {A phase retrieval algorithm for shifting illumination},
  volume = {85},
  ISSN = {1077-3118},
  url = {http://dx.doi.org/10.1063/1.1823034},
  DOI = {10.1063/1.1823034},
  number = {20},
  journal = {Applied Physics Letters},
  publisher = {AIP Publishing},
  author = {Rodenburg,  J. M. and Faulkner,  H. M. L.},
  year = {2004},
  month = Nov,
  pages = {4795--4797}
}

@article{ePIE,
title = {An improved ptychographical phase retrieval algorithm for diffractive imaging},
journal = {Ultramicroscopy},
volume = {109},
number = {10},
pages = {1256--1262},
year = {2009},
issn = {0304-3991},
doi = {https://doi.org/10.1016/j.ultramic.2009.05.012},
url = {https://www.sciencedirect.com/science/article/pii/S0304399109001284},
author = {Andrew M. Maiden and John M. Rodenburg},
}

@article{rPIE,
author = {Andrew Maiden and Daniel Johnson and Peng Li},
journal = {Optica},
number = {7},
pages = {736--745},
publisher = {Optica Publishing Group},
title = {Further improvements to the ptychographical iterative engine},
volume = {4},
month = {Jul},
year = {2017},
url = {https://opg.optica.org/optica/abstract.cfm?URI=optica-4-7-736},
doi = {10.1364/OPTICA.4.000736},
}

@book{lang1985complex,
  title = {Complex Analysis},
  ISBN = {9781475718713},
  ISSN = {0072-5285},
  url = {http://dx.doi.org/10.1007/978-1-4757-1871-3},
  DOI = {10.1007/978-1-4757-1871-3},
  series = {Graduate Texts in Mathematics},
  volume   = {103},
  publisher = {Springer},
  address   = {New York},
  author = {Lang,  Serge},
  year = {1985}
}

@misc{CRcalculus,
  doi = {10.48550/ARXIV.0906.4835},
  url = {https://arxiv.org/abs/0906.4835},
  author = {Kreutz-Delgado,  Ken},
  title = {The Complex Gradient Operator and the {CR}-Calculus},
  publisher = {arXiv},
  year = {2009},
  copyright = {arXiv.org perpetual,  non-exclusive license}
}

@article{Proximal,
  title = {Proximal Algorithms},
  volume = {1},
  ISSN = {2167-3918},
  url = {http://dx.doi.org/10.1561/2400000003},
  DOI = {10.1561/2400000003},
  number = {3},
  journal = {Foundations and Trends in Optimization},
  publisher = {Emerald},
  author = {Parikh,  Neal and Boyd,  Stephen},
  year = {2014},
  month = Jan,
  pages = {127--239}
}

\appendix
\section{Notational convention}

For the sake of readability in the computational details, we adopt a slightly
informal convention for element-wise division that differs from the main
text. Whenever the numerator and denominator have the same dimensions
(e.g., vectors or matrices in $\mathbb{C}^n$), we write
\[
    \frac{\bm{x}}{\bm{y}}
\]
to denote element-wise division, that is,
\[
    \frac{\bm{x}}{\bm{y}} \defeq (x_i / y_i)_i .
\]
Division by a scalar complex number $c$ is always written explicitly as
$\frac{1}{c}\,\bm{x}$, and never using the potentially ambiguous notation
$\bm{x}/c$. No matrix inversion appears in these derivations, and matrix
inverses are never represented using fraction notation. Thus, this slight
abuse of notation is unambiguous. In the main text, element-wise
division is denoted using the standard symbol~$\oslash$.

\section[Computation for the Misfit- and Surrogate-Derivative]{Computation for complex gradients}
We consider the Euclidean space representations $\mathbb{R}^{2m^2}$ on the vectors in the complex field $\bm{z}_k\in \mathbb{C}^{m^2}$ as follows:
\begin{equation*}
\left\{ [\bm{x}_k^\top, \bm{y}_k^\top]^\top ; \bm{x}_k + \bm{y}_k i = \bm{z}_k \right\}\,.
\end{equation*}
With the projection operator:
\begin{equation*}
    \mathcal{R}_k(\bm{z}_k) = \mathcal{F}^{-1}\left[\sqrt{\bm{d}_k} \odot \exp\left(i \theta(\mathcal{F}(\bm{Q}\odot \bm{z}_k))\right)\right]\,,
\end{equation*}
we define the objective function in the $k$-th scanning region:
\begin{equation*}
\Phi_k(\bm{z}_k) = \frac{1}{2} \left\|\bm{Q} \odot\bm{z}_k - \mathcal{R}_k(\bm{z}_k) \right\|_2^2\,,
\end{equation*}
and find its complex gradient:
\begin{equation*}
\nabla_{\bm{z}_k} \Phi_k = \nabla_{\bm{x}_k} \Phi_k + i\nabla_{\bm{y}_k} \Phi_k\,.
\end{equation*}

Firstly, the Parseval's theorem for the discrete Fourier transform gives:
\begin{equation*}
\|\mathcal{F}(\bm{z}_k)\|_2^2=m^2\|\bm{z}_k\|_2^2\quad\text{and}\quad\|\mathcal{F}^{-1}(\bm{z}_k)\|_2^2=\frac{1}{m^2}\|\bm{z}_k\|_2^2\,.
\end{equation*}
Consequently, for all $\bm{v}\in\mathbb{C}^{m^2}$ we have
\begin{equation*}
\langle\mathcal{F}(\bm{v}),\mathcal{F}(\bm{v})\rangle = \langle \bm{v},\mathcal{F}^\ast\mathcal{F}(\bm{v})\rangle = m^2\langle \bm{v},\bm{v}\rangle\,,
\end{equation*}
which implies 
\begin{equation*}
    \frac{1}{m^2}\mathcal{F}^\ast\mathcal{F}(\bm{v}) = \bm{I} \implies \mathcal{F}^{-1} = \frac{1}{m^2}\mathcal{F}^\ast\,.
\end{equation*}

\subsection[Computation for the Misfit Derivative]{Computation for~Eqn.~\ref{eqn:misfit_derivative}}\label{app:computation_misfit_derivative}
We find
\begin{equation*}
\nabla_{\bm{z}_k} \Phi_k = \overline{\bm{Q}}\odot(\bm{Q}\odot\bm{z}_k - \mathcal{R}_k(\bm{z}_k))\,.
\end{equation*}
For the $k$-th scanning region, we have
\begin{equation*}
    \begin{aligned}
        \Phi_k(\bm{z}_k) &= \frac{1}{2} \left\|\bm{Q}\odot \bm{z}_k - \mathcal{R}_k(\bm{z}_k) \right\|_2^2\,,\\
        &= \frac{1}{2m^2} \left\| \mathcal{F}(\bm{Q}\odot\bm{z}_k)  -\sqrt{\bm{d}_k}\odot\exp\left(i \theta(\mathcal{F}(\bm{Q}\odot \bm{z}_k))\right)\right\|_2^2\,,\\
        &= \frac{1}{2m^2} \left\||\mathcal{F}(\bm{Q}\odot\bm{z}_k)| -  \sqrt{\bm{d}_k}\right\|_2^2\,.
    \end{aligned}
\end{equation*}

Using the chain rule for the cogradient operator defined in Chapter 4.2 of~\cite{CRcalculus}, we have
\begin{equation*}
    \begin{aligned}
        \frac{\partial\,  |\mathcal{F}(\bm{Q}\odot\bm{z}_k)|}{\partial\,  \bm{z}_k}&= \frac{\partial\,  |\mathcal{F}(\bm{Q}\odot\bm{z}_k)|}{\partial\,  \mathcal{F}(\bm{Q}\odot\bm{z}_k)}\frac{\partial\,  \mathcal{F}(\bm{Q}\odot\bm{z}_k)}{\partial\,  \bm{z}_k} +  \frac{\partial\,  |\mathcal{F}(\bm{Q}\odot\bm{z}_k)|}{\partial\,  \overline{\mathcal{F}(\bm{Q}\odot\bm{z}_k)}}\frac{\partial\,  \overline{\mathcal{F}(\bm{Q}\odot\bm{z}_k)}}{\partial\,  \bm{z}_k}\,,\\
        &= \operatorname{diag}\left(\frac{\overline{\mathcal{F}(\bm{Q}\odot\bm{z}_k)}}{2|\mathcal{F}(\bm{Q}\odot\bm{z}_k)|}\right)\mathcal{F}\operatorname{diag}(\bm{Q})\,.
    \end{aligned}
\end{equation*}
We note that $\frac{\partial\,  \overline{\mathcal{F}(\bm{Q}\odot\bm{z}_k)}}{\partial\,  \bm{z}_k}$ can be computed as 0.

Consequently, we have
\begin{equation*}
    \begin{aligned}
        \frac{\partial\,  \Phi_{k}}{\partial\,  \bm{z}_k}  &= \frac{1}{2m^2}\frac{\partial\,  \left\||\mathcal{F}(\bm{Q}\odot\bm{z}_k)| -  \sqrt{\bm{d}_k}\right\|_2^2}{\partial\,  \bm{z}_k} \,,\\
        &= \frac{1}{m^2}\left(|\mathcal{F}(\bm{Q}\odot\bm{z}_k)|-\sqrt{\bm{d}_k}\right)\frac{\partial\,  |\mathcal{F}(\bm{Q}\odot\bm{z}_k)|}{\partial\,  \bm{z}_k}\,,\\ 
        &= \frac{1}{m^2}\left(|\mathcal{F}(\bm{Q}\odot\bm{z}_k)|-\sqrt{\bm{d}_k}\right)\operatorname{diag}\left(\frac{\overline{\mathcal{F}(\bm{Q}\odot\bm{z}_k)}}{2|\mathcal{F}(\bm{Q}\odot\bm{z}_k)|}\right)\mathcal{F}\operatorname{diag}(\bm{Q})\,.\\  
    \end{aligned}
\end{equation*}

Finally, by definition, we obtain
\begin{equation*}
\begin{aligned}
    \nabla_{\bm{z}_k} \Phi_k &= \nabla_{\bm{x}_k} \Phi_k + i\nabla_{\bm{y}_k} \Phi_k\,,\\
    &= \left(\frac{\partial\,  \Phi_{k}}{\partial\,  \bm{x}_k} + i\frac{\partial\,  \Phi_{k}}{\partial\,  \bm{y}_k}\right)^\top\,,\\
     &= 2\frac{\partial\,  \Phi_{k}}{\partial\,  \bm{z}_k}^*\,,\\
    &= \operatorname{diag}\left(\overline{\bm{Q}}\right)\left(\frac{1}{m^2}\mathcal{F}^*\right)\left[\operatorname{diag}\left(\frac{\mathcal{F}(\bm{Q}\odot\bm{z}_k)}{|\mathcal{F}(\bm{Q}\odot\bm{z}_k)|}\right) \left(|\mathcal{F}(\bm{Q}\odot\bm{z}_k)|-\sqrt{\bm{d}_k}\right)\right]\,,\\
    &= \overline{\bm{Q}}\odot\mathcal{F}^{-1}\left[\exp\left(i \theta(\mathcal{F}(\bm{Q}\odot \bm{z}_k))\right)\odot\left(|\mathcal{F}(\bm{Q}\odot\bm{z}_k)|-\sqrt{\bm{d}_k}\right)\right]\,,\\
    &= \overline{\bm{Q}}\odot\mathcal{F}^{-1}\left[\mathcal{F}(\bm{Q}\odot\bm{z}_k)-\sqrt{\bm{d}_k}\odot\exp\left(i \theta(\mathcal{F}(\bm{Q}\odot \bm{z}_k))\right)\right]\,,\\
    &= \overline{\bm{Q}}\odot\left[\bm{Q}\odot\bm{z}_k-\mathcal{F}^{-1}\left(\sqrt{\bm{d}_k}\odot\exp\left(i \theta(\mathcal{F}(\bm{Q}\odot \bm{z}_k))\right)\right)\right]\,,\\
    &= \overline{\bm{Q}}\odot\left(\bm{Q}\odot\bm{z}_k-\mathcal{R}_k(\bm{z}_k)\right)\,.
\end{aligned}
\end{equation*}
Note that this computation does not work when $\mathcal{F}(\bm{Q}\odot \bm{z}_k)_r=0$ at some index $r$. This issue can be prevented by finding the gradient using the third statement of Proposition~\ref{prop:mm_cond}. 

\subsection[Computation for the Surrogate Derivative]{Computation for Eqn.~\ref{eqn:derivative_surrogate}}\label{app:computation_surrogate_derivative}
We find 
\begin{equation*} 
    \nabla_{\bm{z}_k}\Phi^{\texttt{PIE}}_k\left(\bm{z}_k ; \bm{z}_k^{(j)}\right) = \overline{\bm{Q}}\odot\left(\bm{Q}\odot\bm{z}_k - \mathcal{R}_k\left(\bm{z}_k^{(j)}\right)\right) + \bm{u}\odot\left(\bm{z}_k - \bm{z}_k^{(j)}\right)\,,
\end{equation*}
where 
\begin{equation*}
    \Phi^{\texttt{PIE}}_k\left(\bm{z}_k ; \bm{z}_k^{(j)}\right) = \frac{1}{2} \left\|\bm{Q}\odot \bm{z}_k - \mathcal{R}_k\left(\bm{z}_k^{(j)}\right)\right\|_2^2 + \frac{1}{2}\bm{u}^\top\left|\bm{z}_k - \bm{z}_k^{(j)}\right|^2\,.
\end{equation*}

For the first term, we have
\begin{equation*}
\begin{aligned}
     \frac{\partial\, \left\|\bm{Q}\odot \bm{z}_k -\mathcal{R}_k\left(\bm{z}_k^{(j)}\right)\right\|_2^2}{\partial\,   \bm{z}_k} &= \frac{\partial\,  \left(\bm{Q}\odot \bm{z}_k -\mathcal{R}_k\left(\bm{z}_k^{(j)}\right)\right)\cdot\overline{\left(\bm{Q}\odot \bm{z}_k -\mathcal{R}_k\left(\bm{z}_k^{(j)}\right)\right)}}{\partial\,   \bm{z}_k}\\
     &= \overline{\left(\bm{Q}\odot \bm{z}_k -\mathcal{R}_k\left(\bm{z}_k^{(j)}\right)\right)}\operatorname{diag}(\bm{Q})\,.
\end{aligned}
\end{equation*}
For the second term, 
\begin{equation*}
\begin{aligned}
     \frac{\partial\,  \bm{u}\cdot\left|\bm{z}_k-\bm{z}_k^{(j)}\right|^2}{\partial\,   \bm{z}_k} &= \frac{\partial\,  \left(\bm{z}_k-\bm{z}_k^{(j)}\right)\cdot\left(\bm{u}\odot\overline{\left(\bm{z}_k-\bm{z}_k^{(j)}\right)}\right)}{\partial\,   \bm{z}_k}\\
     &= \bm{u}\odot\overline{\left(\bm{z}_k-\bm{z}_k^{(j)}\right)}\,.
\end{aligned}
\end{equation*}

Therefore, 
\begin{equation*}
    \begin{aligned}
    \nabla_{\bm{z}_k} \Phi^{\texttt{PIE}}_k &= 2\frac{\partial\,   \Phi^{\texttt{PIE}}_k}{\partial\,  \bm{z}_k}^*\,,\\
    &= 2\cdot\frac{1}{2}\cdot\left(\operatorname{diag}\left(\overline{\bm{Q}}\right)\left(\bm{Q}\odot \bm{z}_k -\mathcal{R}_k\left(\bm{z}_k^{(j)}\right)\right)+\bm{u}\odot\left(\bm{z}_k-\bm{z}_k^{(j)}\right)\right)\,,\\
    &= \overline{\bm{Q}}\odot\left(\bm{Q}\odot \bm{z}_k -\mathcal{R}_k\left(\bm{z}_k^{(j)}\right)\right)+\bm{u}\odot\left(\bm{z}_k-\bm{z}_k^{(j)}\right)\,.
    \end{aligned}
\end{equation*}

\section{Proof for Proposition~\ref{prop:mm_cond} and Proposition~\ref{thm:convergence_iterative_surrogate}}
We begin with a lemma and a proposition.
\begin{lemma}\label{lem:convergence_iterative_surrogate}
    For any complex numbers $x,y\in\mathbb{C}$ and $d\geq0$, 
    \begin{equation*}
        \left|x-de^{i\theta(x)}\right|^2\leq\left|x-de^{  i\theta(y)}\right|^2\,.
    \end{equation*}
    Additionally, equality implies either $d = 0$, $x = 0$, or $\theta(x)=\theta(y)$.
    \begin{proof}
        The first term can be expanded as:
        \begin{equation*}
            \begin{aligned}
                \left|x-de^{i\theta(x)}\right|^2&=\left||x|e^{i\theta(x)}-de^{i\theta(x)}\right|^2\,,\\
                &= \left||x|-d\right|^2\left|e^{i\theta(x)}\right|^2\,,\\
                &=|x|^2+d^2-2|x|d\,.
            \end{aligned}
        \end{equation*}
        The second term can be expanded as:
        \begin{equation*}
            \begin{aligned}
                \left|x-de^{i\theta(y)}\right|^2&=\left||x|e^{i\theta(x)}-de^{i\theta(y)}\right|^2\,,\\
                &=\left(|x|e^{i\theta(x)}-de^{i\theta(y)}\right)\left(|x|e^{-i\theta(x)}-de^{-i\theta(y)}\right)\,,\\
                &=|x|^2+d^2-|x|d\left(e^{i(\theta(x)-\theta(y))}+e^{i(\theta(y)-\theta(x))}\right)\,,\\
                &=|x|^2+d^2-2|x|d\cos(\theta(x)-\theta(y))\\
            \end{aligned}
        \end{equation*}
        The inequality follows from the fact that $\cos(\theta(x)-\theta(y))\leq1$ and that $d\geq0$. Additionally, equality implies either $d=0$, $|x|=0$, or $\cos(\theta(x)-\theta(y))=1$, which is equivalent to either $d=0$, $x=0$, or $\theta(x)=\theta(y)$.
    \end{proof}
\end{lemma}
\begin{proposition}\label{prop:convergence_iterative_surrogate}
For each scanning region indexed by $k=1,\dots,N$,
\begin{equation*}
    \left\|\bm{Q}\odot\bm{z}_k - \mathcal{R}_k\left(\bm{z}_k\right)\right\|_2^2\leq \left\|\bm{Q}\odot\bm{z}_k - \mathcal{R}_k\left(\bm{z}^{(j)}_k\right)\right\|_2^2\,.
\end{equation*}
Additionally, equality implies $\theta\left(\mathcal{F}\left(\bm{Q}\odot \bm{z}_k\right)\right) = \theta\left(\mathcal{F}\left(\bm{Q}\odot \bm{z}_k^{(j)}\right)\right)$ except at the zero entries of $\bm{d}_k$ or at the zero entries of $\mathcal{F}\left(\bm{Q}\odot \bm{z}_k\right)$.

\begin{proof}
    Using Parseval’s theorem for the discrete Fourier transform, it suffices to show:
    \begin{equation*}
        \left\|\mathcal{F}\left(\bm{Q}\odot \bm{z}_k\right) - \sqrt{\bm{d}_k} \odot\exp\left(i \theta\left(\mathcal{F}\left(\bm{Q}\odot \bm{z}_k\right)\right)\right)\right\|_2^2\leq \left\|\mathcal{F}\left(\bm{Q}\odot \bm{z}_k\right) - \sqrt{\bm{d}_k}\odot\exp\left(i \theta\left(\mathcal{F}\left(\bm{Q}\odot \bm{z}_k^{(j)}\right)\right)\right)\right\|_2^2\,.
    \end{equation*}
    The inequality and the additional equality outcomes follow from Lemma~\ref{lem:convergence_iterative_surrogate}, which applies to the inequality element-wise. 
\end{proof}
\end{proposition}

\subsection{Proof for Proposition~\ref{prop:mm_cond}}\label{app:proof_mm_cond}
\begin{proof}
     It is straightforward to check the first statement.  
     
     For the second statement, it follows from Proposition~\ref{prop:convergence_iterative_surrogate} that  $\Phi_k\left(\bm{z}_k\right) \leq \widetilde{\Phi}_k\left(\bm{z}_k ; \bm{z}_k^{(j)}\right)$ for each scanning region.  Hence,
     \begin{equation*}
         \Phi\left(\bm{z}\right)=\sum_{k=1}^N\Phi_k\left(\bm{z}_k\right) \leq \sum_{k=1}^N\widetilde{\Phi}_k\left(\bm{z}_k ; \bm{z}_k^{(j)}\right)=\widetilde{\Phi}\left(\bm{z} ; \bm{z}^{(j)}\right)\,.
     \end{equation*}

     For the last statement, we consider equivalent representations of the two objective functions using the real and imaginary parts of the input $\bm{z}=\bm{x}+i\bm{y}$: 
     \begin{equation*}
         \Phi\left(\bm{z}\right) \equiv \Phi\left(\bm{x},\bm{y}\right)\quad\text{and}\quad \widetilde{\Phi}\left(\bm{z};\bm{z}^{(j)}\right) \equiv \widetilde{\Phi}\left(\bm{x},\bm{y};\bm{z}^{(j)}\right)\,.
     \end{equation*}
    With $\bm{z}^{(j)}=\bm{x}^{(j)}+i\bm{y}^{(j)}$, it suffices to show:
    \begin{equation*}
        \nabla_{\bm{x}}\Phi\left(\bm{x}^{(j)},\bm{y}^{(j)}\right) = \nabla_{\bm{x}}\widetilde{\Phi}\left(\bm{x}^{(j)},\bm{y}^{(j)};\bm{z}^{(j)}\right)\quad\text{and}\quad \nabla_{\bm{y}}\Phi\left(\bm{x}^{(j)},\bm{y}^{(j)}\right) = \nabla_{\bm{y}}\widetilde{\Phi}\left(\bm{x}^{(j)},\bm{y}^{(j)};\bm{z}^{(j)}\right)\,.
    \end{equation*}
    By the previous two statements and definition of the gradient, for a small perturbation $\epsilon \bm{h}$ such that $\epsilon > 0$ and $\bm{h}\in\mathbb{R}^{m^2}$, we have
    \begin{equation*}
        \begin{aligned}
            \epsilon\left\langle\nabla_{\bm{x}}\Phi\left(\bm{x}^{(j)},\bm{y}^{(j)}\right),\bm{h}\right\rangle+o(\epsilon)&= \Phi\left(\bm{x}^{(j)}+\epsilon \bm{h},\bm{y}^{(j)}\right)-\Phi\left(\bm{x}^{(j)},\bm{y}^{(j)}\right)\,,\\
            &\leq \widetilde{\Phi}\left(\bm{x}^{(j)}+\epsilon\bm{h},\bm{y}^{(j)};\bm{z}^{(j)}\right)-\widetilde{\Phi}\left(\bm{x}^{(j)},\bm{y}^{(j)};\bm{z}^{(j)}\right)\\
            &=\epsilon\left\langle\nabla_{\bm{x}}\widetilde{\Phi}\left(\bm{x}^{(j)},\bm{y}^{(j)};\bm{z}^{(j)}\right),\bm{h}\right\rangle+o(\epsilon)\,.
        \end{aligned}
    \end{equation*}
    Noting that $\epsilon > 0$ and $\bm{h}\in\mathbb{R}^{m^2}$ are arbitrary, we have
    \begin{equation*}
        \nabla_{\bm{x}}\Phi\left(\bm{x}^{(j)},\bm{y}^{(j)}\right) = \nabla_{\bm{x}}\widetilde{\Phi}\left(\bm{x}^{(j)},\bm{y}^{(j)};\bm{z}^{(j)}\right)\,.
    \end{equation*}
    Similarly, we can show the equality for $\nabla_{\bm{y}}$ and conclude that 
    \begin{equation*}
        \nabla_{\bm{z}}\Phi\left(\bm{z}^{(j)}\right) = \nabla_{\bm{z}}\widetilde{\Phi}\left(\bm{z}^{(j)};\bm{z}^{(j)}\right)\,.
    \end{equation*}
\end{proof}

\subsection{Proof for Proposition~\ref{thm:convergence_iterative_surrogate}}\label{app:proof_convergence_iterative_surrogate}

\begin{proof}
We separate the proof into two parts:
\begin{enumerate}
    \item \textbf{Inequality}: Since we obtain $\bm{z}^{(j+1)}$ by solving
    \begin{equation*}
        \bm{z}^{(j+1)} = \operatornamewithlimits{argmin}_{\bm{z}'} \frac{1}{2} \sum_{k=1}^N \left\|\bm{Q}\odot\bm{z}'_k - \mathcal{R}_k\left(\bm{z}^{(j)}_k\right)\right\|_2^2\,,
    \end{equation*}
    we have
    \begin{equation}\label{eqn:proof_thm_1}
        \frac{1}{2}\sum_{k=1}^N \left\|\bm{Q}\odot\bm{z}^{(j+1)}_k - \mathcal{R}_k\left(\bm{z}^{(j)}_k\right)\right\|_2^2 \leq \frac{1}{2}\sum_{k=1}^N \left\|\bm{Q}\odot\bm{z}^{(j)}_k - \mathcal{R}_k\left(\bm{z}^{(j)}_k\right)\right\|_2^2=\Phi\left(\bm{z}^{(j)}\right)\,.
    \end{equation}

    Then, by Proposition~\ref{prop:convergence_iterative_surrogate}, we have
    \begin{equation}\label{eqn:proof_thm_2}
        \Phi\left(\bm{z}^{(j+1)}\right)=\frac{1}{2}\sum_{k=1}^N \left\|\bm{Q}\odot\bm{z}^{(j+1)}_k - \mathcal{R}_k\left(\bm{z}^{(j+1)}_k\right)\right\|_2^2\leq\frac{1}{2}\sum_{k=1}^N \left\|\bm{Q}\odot\bm{z}^{(j+1)}_k - \mathcal{R}_k\left(\bm{z}^{(j)}_k\right)\right\|_2^2\,.
    \end{equation}
Combining Eqn.~\eqref{eqn:proof_thm_1} and Eqn.~\eqref{eqn:proof_thm_2} completes the proof of the inequality.
\item \textbf{Optimality:} By Proposition~\ref{prop:convergence_iterative_surrogate}, the inequality~Eqn.~\eqref{eqn:proof_thm_2} achieves equality only when $\theta\left(\mathcal{F}\left(\bm{Q}\odot \bm{z}_k^{(j+1)}\right)\right) = \theta\left(\mathcal{F}\left(\bm{Q}\odot \bm{z}_k^{(j)}\right)\right)$ except at the zero entries of $\bm{d}_k$ or at the zero entries of $\mathcal{F}\left(\bm{Q}\odot \bm{z}_k^{(j+1)}\right)$ for all $k=1,2,\dots,N$.

Since at the zero entries $\mathcal{F}\left(\bm{Q}\odot \bm{z}_k^{(j+1)}\right)$, we impose $\theta\left(\mathcal{F}\left(\bm{Q}\odot \bm{z}_k^{(j+1)}\right)\right) = \theta\left(\mathcal{F}\left(\bm{Q}\odot \bm{z}_k^{(j)}\right)\right)$ in Eqn.~\ref{eqn:specification}, we have, for all $k=1,2,\dots,N$,  
\begin{equation*}
    \begin{aligned}
        \mathcal{R}_k\left(\bm{z}^{(j+1)}_k\right) &= \mathcal{F}^{-1}\left(\sqrt{\bm{d}_k}\odot \exp\left(i \theta\left(\mathcal{F}\left(\bm{Q}\odot \bm{z}_k^{(j+1)}\right)\right)\right)\right)\,,\\
        &= \mathcal{F}^{-1}\left(\sqrt{\bm{d}_k}\odot \exp\left(i \theta\left(\mathcal{F}\left(\bm{Q}\odot \bm{z}_k^{(j)}\right)\right)\right)\right)\,,\\
        &= \mathcal{R}_k\left(\bm{z}^{(j)}_k\right)\,.
    \end{aligned}
\end{equation*}

That leads to
\begin{equation*}
\begin{aligned}
    \bm{z}^{(j+1)}
    &=\operatornamewithlimits{argmin}_{\bm{z}'} \frac{1}{2} \sum_{k=1}^N \left\|\bm{Q}\odot\bm{z}'_k - \mathcal{R}_k\left(\bm{z}^{(j)}_k\right)\right\|_2^2\,,\\
    &=\operatornamewithlimits{argmin}_{\bm{z}'} \frac{1}{2} \sum_{k=1}^N \left\|\bm{Q}\odot\bm{z}'_k - \mathcal{R}_k\left(\bm{z}^{(j+1)}_k\right)\right\|_2^2\,,\\
     &=\operatornamewithlimits{argmin}_{\bm{z}'} \widetilde{\Phi}\left(\bm{z}';\bm{z}^{(j+1)}\right)\,.\\
\end{aligned}
\end{equation*}

By the optimality of a quadratic objective function, we have
\begin{equation*}
    \nabla_{\bm{z}} \widetilde{\Phi}\left(\bm{z}^{(j+1)};\bm{z}^{(j+1)}\right) = \overline{\bm{Q}}\odot\left(\bm{Q}\odot \bm{z}_k^{(j+1)} -\mathcal{R}_k\left(\bm{z}_k^{(j+1)}\right)\right) = 0\,.
\end{equation*}

Therefore, we can conclude that
\begin{equation*}
    \nabla_{\bm{z}}\Phi\left(\bm{z}^{(j+1)}\right) = \nabla_{\bm{z}} \widetilde{\Phi}\left(\bm{z}^{(j+1)};\bm{z}^{(j+1)}\right) = 0\,.
\end{equation*}

\end{enumerate}
\end{proof}

\section{Proof for Proposition~\ref{thm:convergence_rate}}\label{app:proof_convergence_rate}
\begin{proof}
    By Proposition~\ref{prop:mm_cond} and Eqn.~\ref{eqn:surrogate_minimization}, and Hölder's inequality, we have
    \begin{equation*}
    \begin{aligned}
        \Phi\left(\bm{z}^{(j+1)}\right)-\Phi\left(\bm{z}^{(j)}\right)&\leq \widetilde{\Phi}\left(\bm{z}^{(j+1)} ; \bm{z}^{(j)}\right) - \widetilde{\Phi}\left(\bm{z}^{(j)} ; \bm{z}^{(j)}\right)\,,\\
        &\leq \widetilde{\Phi}\left(\bm{z} ; \bm{z}^{(j)}\right) - \widetilde{\Phi}\left(\bm{z}^{(j)} ; \bm{z}^{(j)}\right)\,,\\
        &= \frac{1}{2} \sum_{k=1}^N \left[\left\|\bm{Q}\odot P_k\bm{z} - \mathcal{R}_k\left(P_k\bm{z}^{(j)}\right)\right\|_2^2-\left\|\bm{Q}\odot P_k\bm{z}^{(j)} - \mathcal{R}_k\left(P_k\bm{z}^{(j)}\right)\right\|_2^2\right]\,,\\
        &= \frac{1}{2} \sum_{k=1}^N \bigg[\left\|\bm{Q}\odot P_k\bm{z}\right\|^2  - \overline{\bm{Q}\odot P_k\bm{z}}\cdot \mathcal{R}_k\left(P_k\bm{z}^{(j)}\right) -\bm{Q}\odot P_k\bm{z}\cdot \overline{\mathcal{R}_k\left(P_k\bm{z}^{(j)}\right)}\\
        &\quad\quad\quad-\left\|\bm{Q}\odot P_k\bm{z}^{(j)}\right\|^2+ \overline{\bm{Q}\odot P_k\bm{z}^{(j)}}\cdot \mathcal{R}_k\left(P_k\bm{z}^{(j)}\right) + \bm{Q}\odot P_k\bm{z}^{(j)}\cdot \overline{\mathcal{R}_k\left(P_k\bm{z}^{(j)}\right)}\bigg]\,,\\
        &= \frac{1}{2} \sum_{k=1}^N \bigg[\left\|\bm{Q}\odot P_k\bm{z}-\bm{Q}\odot P_k\bm{z}^{(j)}\right\|^2 +  \overline{\bm{Q}\odot P_k\bm{z}-\bm{Q}\odot P_k\bm{z}^{(j)}}\cdot \left(\bm{Q}\odot P_k\bm{z}^{(j)} - \mathcal{R}_k\left(P_k\bm{z}^{(j)}\right)\right)\\
        &\quad\quad\quad+\left(\bm{Q}\odot P_k\bm{z}-\bm{Q}\odot P_k\bm{z}^{(j)}\right)\cdot \overline{\bm{Q}\odot P_k\bm{z}^{(j)} - \mathcal{R}_k\left(P_k\bm{z}^{(j)}\right)}\bigg]\,,\\
        &= \frac{1}{2} \sum_{k=1}^N  \bigg[\left(P_k^\top |\bm{Q}|^2\right)\cdot \left|\bm{z}-\bm{z}^{(j)}\right|^2  +  \overline{\bm{z}-\bm{z}^{(j)}}\cdot P_k^{\top}\nabla_{\bm{z}_k}\Phi_k\left(\bm{z}_k^{(j)}\right)\\
        &\quad\quad\quad\quad+\left(\bm{z}-\bm{z}^{(j)}\right)\cdot \overline{P_k^{\top}\nabla_{\bm{z}_k}\Phi_k\left(\bm{z}_k^{(j)}\right)}\bigg]\,,\\
        &= \left(\frac{1}{2} \sum_{k=1}^N  P_k^\top|\bm{Q}|^2\right)\cdot \left|\bm{z}-\bm{z}^{(j)}\right|^2 + \overline{\bm{z}-\bm{z}^{(j)}}\cdot\left(\frac{1}{2} \sum_{k=1}^N P_k^{\top}\nabla_{\bm{z}_k}\Phi_k\left(\bm{z}_k^{(j)}\right)\right)\\
        &\quad\quad\quad+\left(\bm{z}-\bm{z}^{(j)}\right)\cdot\overline{\frac{1}{2} \sum_{k=1}^N P_k^{\top}\nabla_{\bm{z}_k}\Phi_k\left(\bm{z}_k^{(j)}\right)} \,,\\
        &\leq \frac{1}{2}\left\|\sum_{k=1}^N  P_k^\top|\bm{Q}|^2\right\|_\infty\left\|\bm{z}-\bm{z}^{(j)}\right\|_2^2+\frac{1}{2}\left(\overline{\bm{z}-\bm{z}^{(j)}}\right)\cdot \nabla_{\bm{z}}\Phi\left(\bm{z}^{(j)}\right)\\
        &\quad\quad\quad+\frac{1}{2}\left(\bm{z}-\bm{z}^{(j)}\right)\cdot\overline{\nabla_{\bm{z}}\Phi\left(\bm{z}^{(j)}\right)}\,.
    \end{aligned} 
    \end{equation*}
    We denote $\frac{1}{2}\left\|\sum_{k=1}^N  P_k^\top|\bm{Q}|^2\right\|_\infty$ by $L$. Then, the choice $\bm{z}=\bm{z}^{(j)}-\frac{1}{2L}\nabla_{\bm{z}}\Phi\left(\bm{z}^{(j)}\right)$ yields: 
    \begin{equation*}
        \Phi\left(\bm{z}^{(j)}\right)-\Phi\left(\bm{z}^{(j+1)}\right)\geq\frac{1}{4L}\left\|\nabla_{\bm{z}}\Phi\left(\bm{z}^{(j)}\right)\right\|_2^2\,.
    \end{equation*}
    A telescoping sum gives:
    \begin{equation*}
    \begin{aligned}
        \frac{n+1}{4L}\min_{0\leq j\leq n} \left\|\nabla_{\bm{z}}\Phi\left(\bm{z}^{(j)}\right)\right\|_2^2&\leq \frac{1}{4L}\sum_{j=0}^n \left\|\nabla_{\bm{z}}\Phi\left(\bm{z}^{(j)}\right)\right\|_2^2\,,\\
        &\leq \Phi\left(\bm{z}^{(0)}\right)-\Phi\left(\bm{z}^{(n+1)}\right)\,,\\
        &\leq \Phi\left(\bm{z}^{(0)}\right)\,.
    \end{aligned}
    \end{equation*}
\end{proof}

\section{Useful Lemmas} 
\begin{lemma}\label{lem:commutative}
For any $\bm{A}\in\mathbb{C}^{n_h^2}$ and $\bm{B}\in\mathbb{C}^{n_H^2}$, we have
    \begin{equation*}
        \bm{B}\odot \bm{I}_h^H\bm{A} = \bm{I}_h^H((\bm{I}_H^h\bm{B})\odot\bm{A})\,.
    \end{equation*}
Note that when $\bm{A}=\bm{J}$, a vector of ones, we have 
\begin{equation*}
    \bm{B} = \bm{I}_h^H\bm{I}_H^h\bm{B}\implies \bm{I}_h^H\bm{I}_H^h = \bm{I},\text{ the identity matrix}\,.
\end{equation*}

    \begin{proof}
    Note that using the binning strategy, each column of $\bm{I}_h^H$ is non-zero only at one entry, i.e.,
  \begin{equation*}
    (\bm{I}_h^H)_{ij}(\bm{I}_h^H)_{kj} = \frac{1}{4}(\bm{I}_h^H)_{ij}\delta_{ik}\,,
  \end{equation*}
  and by definition, we have
  \begin{equation*}
      \bm{I}_H^h = 4(\bm{I}_h^H)^\top\,.
  \end{equation*}
    For all entries in the $\texttt{LHS}$ and $\texttt{RHS}$, we have 
    \begin{equation*}
            \texttt{LHS}_i = \bm{B}_i(\bm{I}_h^H\bm{A})_i = \bm{B}_i\sum_{j}(\bm{I}_h^H)_{ij}\bm{A}_j\,.
    \end{equation*}
    \item 
    \begin{equation*}
        \begin{aligned}
            \texttt{RHS}_i &= (\bm{I}_h^H[(\bm{I}_H^h\bm{B})\odot\bm{A}])_i\,,\\
            &= \sum_{j}(\bm{I}_h^H)_{ij}[(\bm{I}_H^h\bm{B})\odot\bm{A}]_j\,,\\
            &= \sum_{j}(\bm{I}_h^H)_{ij}\left[\sum_{k}(\bm{I}_H^h)_{jk}\bm{B}_k\right]\bm{A}_j\,,\\
            &= \sum_{k}\bm{B}_k\sum_{j}(\bm{I}_h^H)_{ij}(\bm{I}_H^h)_{jk}\bm{A}_j\,,\\
            &= \sum_{k}\bm{B}_k\sum_{j}(\bm{I}_h^H)_{ij}4(\bm{I}_h^H)_{kj}\bm{A}_j\,,\\
             &= \sum_{k}\bm{B}_k\sum_{j}(\bm{I}_h^H)_{ij}\delta_{ik}\bm{A}_j\,,\\
             &= \bm{B}_i\sum_{j}(\bm{I}_h^H)_{ij}\bm{A}_j\,.
        \end{aligned}
    \end{equation*}
    \end{proof}
\end{lemma}

\begin{lemma}\label{lem:unary}
    For any function $f:\mathbb{C}\to\mathbb{C}$, we define the operators
\begin{equation*}
    \bm{f}_H:\mathbb{C}^{n_H^2}\to\mathbb{C}^{n_H^2}\quad\text{and}\quad \bm{f}_h:\mathbb{C}^{n_h^2}\to\mathbb{C}^{n_h^2}\,,
\end{equation*}
to be the element-wise applications of $f$, i.e., for any input vector, the operator applies $f$ to every entry. Then, for any $\bm{B}\in\mathbb{C}^{n_H^2}$,
\begin{equation*}
    \bm{f}_h\left(\bm{I}_H^h\bm{B}\right)= \bm{I}_H^h\bm{f}_H(\bm{B})\,.
\end{equation*}
\end{lemma}

\begin{lemma}\label{lem:binary}
    For any $\bm{B}$, and $\bm{C}\in\mathbb{C}^{n_H^2}$ and any binary function $g:\mathbb{C}\otimes\mathbb{C}\to\mathbb{C}$, we define the operators
\begin{equation*}
    \bm{g}_H:\mathbb{C}^{n_H^2}\otimes\mathbb{C}^{n_H^2}\to\mathbb{C}^{n_H^2}\quad\text{and}\quad \bm{g}_h:\mathbb{C}^{n_h^2}\otimes\mathbb{C}^{n_h^2}\to\mathbb{C}^{n_h^2}\,,
\end{equation*}
to be the element-wise applications of $g$. Then, for any $\bm{B},\bm{C}\in\mathbb{C}^{n_H^2}$,
\begin{equation*}
    \bm{g}_h\left(\bm{I}_H^h\bm{B},\bm{I}_H^h\bm{C}\right)= \bm{I}_H^h\bm{g}_H(\bm{B},\bm{C})\,.
\end{equation*}
\end{lemma}

Lemma~\ref{lem:unary} and Lemma~\ref{lem:binary} are straightforward to check by comparing corresponding entries of both sides, so we omit the proofs.

\begin{lemma}\label{lem:operator_norm}
The operator norms of $\bm{I}_h^H$ and $\bm{I}_H^h$ induced by the 2-norm are:
\begin{equation*}
   \left\|\bm{I}_h^H\right\|_2=\frac{1}{2}\quad\text{and}\quad\left\|\bm{I}_H^h\right\|_2=2\,.
\end{equation*}

The operator norms of $\bm{I}_h^H$ and $\bm{I}_H^h$ induced by the $\infty$-norm are:
\begin{equation*}
   \left\|\bm{I}_h^H\right\|_\infty=1\quad\text{and}\quad\left\|\bm{I}_H^h\right\|_\infty=1\,.
\end{equation*}
    \begin{proof}
       For any vector $\bm{v}\in \mathbb{C}^{n_H^2}$, the prolongation operator $\bm{I}_H^h$ interpolates $\bm{v}$ by replicating each entry uniformly across its corresponding bin. Hence, 
    \begin{equation*}
        \left\|\bm{I}_H^h\bm{v}\right\|_2=2\left\|\bm{v}\right\|_2\quad\text{and}\quad\left\|\bm{I}_H^h\bm{v}\right\|_\infty=\left\|\bm{v}\right\|_\infty\,.
    \end{equation*}
       This leads to 
    \begin{equation*}
        \left\|\bm{I}_H^h\right\|_2=2\quad\text{and}\quad\left\|\bm{I}_H^h\right\|_\infty=1\,.
    \end{equation*}

    For any vector $\bm{w}\in \mathbb{C}^{n_h^2}$, the restriction operator $\bm{I}_h^H$ downsamples $\bm{w}$ by taking average in bins of 4 entries. Let $\mathcal{B}_j$ denote the set of indices of the entries in the bin corresponding to the entry $\left(\bm{I}_h^H\bm{w}\right)_j$.  Then, by Cauchy–Schwarz inequality, 
    \begin{equation*}
         \left|\left(\bm{I}_h^H\bm{w}\right)_j\right| = \left|\frac{1}{4}\sum_{i\in\mathcal{B}_j}\bm{w}_i\right|\leq \frac{1}{2}\sqrt{\sum_{i\in\mathcal{B}_j}|\bm{w}_i|^2}\,.
    \end{equation*}   
    Hence, we have
    \begin{equation*}
        \left\|\bm{I}_h^H\bm{w}\right\|_2^2 = \sum_{j=1}^{n_H^2} \left|\left(\bm{I}_h^H\bm{w}\right)_j\right|^2
        \leq \frac{1}{4}\sum_{j=1}^{n_H^2}\sum_{i\in\mathcal{B}_j}|\bm{w}_i|^2
        = \frac{1}{4}\left\|\bm{w}\right\|_2^2\,.
    \end{equation*}
    Notice that equality is achieved when $\bm{w}$ is a vector of ones, by definition, we have
    \begin{equation*}
        \left\|\bm{I}_h^H\right\|_2=\frac{1}{2}\,.
    \end{equation*}
    
    Similarly, by triangle inequality,
    \begin{equation*}
         \left|\left(\bm{I}_h^H\bm{w}\right)_j\right| = \left|\frac{1}{4}\sum_{i\in\mathcal{B}_j}\bm{w}_i\right|\leq \frac{1}{4}\sum_{i\in\mathcal{B}_j}|\bm{w}_i|\leq\max_{i\in\mathcal{B}_j}|\bm{w}_i|\,,
    \end{equation*}  
    which leads to
    \begin{equation*}
        \left\|\bm{I}_h^H\bm{w}\right\|_\infty = \max_{j=1,\dots,n_H^2}\left|\left(\bm{I}_h^H\bm{w}\right)_j\right|
        \leq \max_{i=1,\dots,n_h^2}|\bm{w}_i|
        = \left\|\bm{w}\right\|_\infty\,.
    \end{equation*}
    The equality is also achieved when $\bm{w}$ is a vector of ones.  Hence, we have
    \begin{equation*}
        \left\|\bm{I}_h^H\right\|_\infty=1\,.
    \end{equation*}
    \end{proof}
\end{lemma}

\section{Proof for Proposition~\ref{prop:weights_property}}\label{app:proof_weights_property}

\begin{proof} We provide proofs for the statements in the same order. 
\begin{itemize}
    \item By Lemma~\ref{lem:commutative} and Lemma~\ref{lem:unary}, we have
    \begin{equation*}
        \begin{aligned}
            \bm{I}_h^H{\bm{W}_{\bm{z}}} &= \bm{I}_h^H\left( \frac{|\bm{Q}|^2}{\bm{I}_H^h\bm{I}_h^H|\bm{Q}|^2}\right)\,,\\
            &= \bm{I}_h^H\left( \bm{I}_H^h\left(\frac{1}{\bm{I}_h^H|\bm{Q}|^2}\right)\odot|\bm{Q}|^2\right)\,,\\
            &= \left(\frac{1}{\bm{I}_h^H|\bm{Q}|^2}\right)\odot\bm{I}_h^H|\bm{Q}|^2\,,\\
            &= \bm{J}\,,\quad\text{a vector of ones}\,.
        \end{aligned}
    \end{equation*}
    For the sake of a contradiction, we assume that $(\bm{W}_{\bm{z}})_r > 4$ at some index $r$. Since $\bm{W}_{\bm{z}}\geq 0$, we have $\bm{I}_h^H{\bm{W}_{\bm{z}}}\neq \bm{J}$ because in the bin of index $r$, the average is greater than 1. That is a contradiction to the first statement. Therefore, $\bm{W}_{\bm{z}}$ element-wise bounded by $4$.
    \item First, we find that
    \begin{equation*}
        \left|\frac{(\bm{I}_H^h\bm{Q}_H)\odot\bm{W}_{\bm{z}}}{\bm{Q}}\right| =  \frac{|(\bm{I}_H^h\bm{Q}_H)|\odot|\bm{Q}|^2}{|\bm{Q}|\odot\left(\bm{I}_H^h\bm{I}_h^H|\bm{Q}|^2\right)} = \frac{\left(\bm{I}_H^h|\bm{Q}_H|\right)\odot|\bm{Q}|}{\bm{I}_H^h\bm{I}_h^H|\bm{Q}|^2}\leq \frac{\left(\bm{I}_H^h\bm{I}_h^H|\bm{Q}|\right)\odot|\bm{Q}|}{\bm{I}_H^h\bm{I}_h^H|\bm{Q}|^2} \,.
    \end{equation*}
    Then, by Lemma~\ref{lem:commutative}, Lemma~\ref{lem:unary}, and Lemma~\ref{lem:binary}, we have
    \begin{equation*}
        \begin{aligned}
            \bm{I}_h^H\left|\frac{(\bm{I}_H^h\bm{Q}_H)\odot\bm{W}_{\bm{z}}}{\bm{Q}}\right|&\leq\bm{I}_h^H\left(\frac{\left(\bm{I}_H^h\bm{I}_h^H|\bm{Q}|\right)\odot|\bm{Q}|}{\bm{I}_H^h\bm{I}_h^H|\bm{Q}|^2}\right)\,,\\&=\bm{I}_h^H\left(\bm{I}_H^h\left(\frac{\bm{I}_h^H|\bm{Q}|}{\bm{I}_h^H|\bm{Q}|^2}\right)\odot|\bm{Q}|\right)\,,\\
            &=\left(\frac{\bm{I}_h^H|\bm{Q}|}{\bm{I}_h^H|\bm{Q}|^2}\right)\odot\bm{I}_h^H|\bm{Q}|\,,\\
            &=\frac{(\bm{I}_h^H|\bm{Q}|)^2}{\bm{I}_h^H|\bm{Q}|^2}\leq 1\,,
        \end{aligned}
    \end{equation*}
    where the last element-wise inequality follows from Jensen's inequality applied to each bin. Hence, by exactly the same argument as the proof of the previous statement, we have
    \begin{equation*}
        \bm{I}_h^H\left|\frac{(\bm{I}_H^h\bm{Q}_H)\odot \bm{W}_{\bm{z}}}{\bm{Q}}\right|\leq 1 \implies \left|\frac{(\bm{I}_H^h\bm{Q}_H)\odot\bm{W}_{\bm{z}}}{\bm{Q}}\right|\leq 4\,,
    \end{equation*}
    where the inequalities are satisfied element-wise. We note the sharp bound for the last inequality is $3/2$, but its more involved derivation is omitted here.
    \item By applying Jensen's inequality to each bin, 
    \begin{equation*}
        |\bm{Q}_H|^2 \leq\left(\bm{I}_h^H|\bm{Q}|\right)^2\leq \bm{I}_h^H |\bm{Q}|^2\implies \frac{|\bm{Q}_H|^2}{\bm{I}_h^H|\bm{Q}|^2}
        \leq 1\,.
    \end{equation*}
\end{itemize}
\end{proof}

\section{Proof for Proposition~\ref{prop:consistency}}\label{app:proof_consistency}
\begin{proof}
By definition, we have
\begin{equation*}
    \begin{aligned}
        \widetilde{\Phi}_{H,k}\left(\bm{I}_h^H\left(\bm{W}_{\bm{z}}\odot\bm{z}_{k}\right); \bm{z}_k^{(j)}\right)
        &= \frac{1}{2} \left\|\bm{Q}_H\odot\bm{I}_h^H\left(\bm{W}_{\bm{z}}\odot\bm{z}_k\right)- \bm{I}_h^H\left(\bm{W}_{\mathcal{R}}\odot\mathcal{R}_k\left(\bm{z}_k^{(j)}\right)\right) \right\|_2^2\,,\\  
        &= \frac{1}{2} \left\|\bm{I}_h^H\left((\bm{I}_H^h\bm{Q}_H)\odot\bm{W}_{\bm{z}}\odot\bm{z}_k\right)- \bm{I}_h^H\left(\bm{W}_{\mathcal{R}}\odot\mathcal{R}_k\left(\bm{z}_k^{(j)}\right)\right) \right\|_2^2\,,\\  
        &=\frac{1}{2} \left\|\bm{I}_h^H\left(\frac{(\bm{I}_H^h\bm{Q}_H)\odot\bm{W}_{\bm{z}}}{\bm{Q}}\odot\bm{Q}\odot \bm{z}_k\right)- \bm{I}_h^H\left(\bm{W}_{\mathcal{R}}\odot\mathcal{R}_k\left(\bm{z}_k^{(j)}\right)\right) \right\|_2^2\,,\\    
        &=\frac{1}{2} \left\|\bm{I}_h^H\left(\bm{W}_{\mathcal{R}}\odot\bm{Q}\odot \bm{z}_k\right)- \bm{I}_h^H\left(\bm{W}_{\mathcal{R}}\odot\mathcal{R}_k\left(\bm{z}_k^{(j)}\right)\right) \right\|_2^2\,,\\ 
        &=\frac{1}{2} \left\|\bm{I}_h^H\left(\bm{W}_{\mathcal{R}}\odot\left(\bm{Q}\odot \bm{z}_k-\mathcal{R}_k\left(\bm{z}_k^{(j)}\right)\right)\right)\right\|_2^2\,,\\ 
        &\leq \frac{1}{2}\left\|\bm{I}_h^H\right\|_2^2\left\|\bm{W}_{\mathcal{R}}\right\|_{\infty}^2\left\|\bm{Q}\odot \bm{z}_k- \mathcal{R}_k\left(\bm{z}_k^{(j)}\right)\right\|_2^2\,,\\
        &=\frac{\left\|\bm{W}_{\mathcal{R}}\right\|_{\infty}^2}{4}\widetilde{\Phi}_k\left(\bm{z}_{k}; \bm{z}_k^{(j)}\right)\,.
    \end{aligned}
\end{equation*}
where the second equality follows from Lemma~\ref{lem:commutative}, and the last equality follows from Lemma~\ref{lem:operator_norm}.

Similarly, we have
\begin{equation*}
    \begin{aligned}
        \left\|\nabla_{\bm{z}_{H,k}}\widetilde{\Phi}_{H,k}\left(\bm{I}_h^H\left(\bm{W}_{\bm{z}}\odot\bm{z}_{k}\right); \bm{z}_k^{(j)}\right)\right\|_2&=\left\|\overline{\bm{Q}_H}\odot\left(\bm{Q}_H\odot\bm{I}_h^H\left(\bm{W}_{\bm{z}}\odot\bm{z}_k\right)- \bm{I}_h^H\left(\bm{W}_{\mathcal{R}}\odot\mathcal{R}_k\left(\bm{z}_k^{(j)}\right)\right)\right) \right\|_2\,,\\
        &=\left\|\overline{\bm{Q}_H}\odot\bm{I}_h^H\left(\bm{W}_{\mathcal{R}}\odot\left(\bm{Q}\odot \bm{z}_k-\mathcal{R}_k\left(\bm{z}_k^{(j)}\right)\right)\right) \right\|_2\,,\\
        &=\left\|\bm{I}_h^H\left(\frac{\left(\bm{I}_H^h\overline{\bm{Q}_H}\right)\odot\bm{W}_{\mathcal{R}}}{\overline{\bm{Q}}}\odot\overline{\bm{Q}}\odot\left(\bm{Q}\odot \bm{z}_k-\mathcal{R}_k\left(\bm{z}_k^{(j)}\right)\right)\right) \right\|_2\,,\\
        &=\left\|\bm{I}_h^H\left(\frac{\left(\bm{I}_H^h\overline{\bm{Q}_H}\right)\odot\bm{W}_{\mathcal{R}}}{\overline{\bm{Q}}}\odot\nabla_{\bm{z}_{k}}\widetilde{\Phi}_k\left(\bm{z}_{k}; \bm{z}_k^{(j)}\right)\right) \right\|_2\,,\\
        &\leq \left\|\bm{I}_h^H\right\|_2\left\|\frac{\left(\bm{I}_H^h\overline{\bm{Q}_H}\right)\odot\bm{W}_{\mathcal{R}}}{\overline{\bm{Q}}}\right\|_{\infty}\left\|\nabla_{\bm{z}_{k}}\widetilde{\Phi}_k\left(\bm{z}_{k}; \bm{z}_k^{(j)}\right) \right\|_2\,.\\
    \end{aligned}
\end{equation*}
Since $\left\|\bm{I}_h^H\right\|_2=\frac{1}{2}$ by Lemma~\ref{lem:operator_norm}, it suffices to check $\left\|\frac{\left(\bm{I}_H^h\overline{\bm{Q}_H}\right)\odot\bm{W}_{\mathcal{R}}}{\overline{\bm{Q}}}\right\|_{\infty} = \left\|\bm{W}_{\bm{u}}^H\right\|_{\infty} $ and it is bounded by 1. To this end, by 
Lemma~\ref{lem:unary}, Lemma~\ref{lem:binary}, we find
\begin{equation*}
    \begin{aligned}
        \left\|\frac{\left(\bm{I}_H^h\overline{\bm{Q}_H}\right)\odot\bm{W}_{\mathcal{R}}}{\overline{\bm{Q}}}\right\|_{\infty}&= \left\|\frac{\overline{\bm{I}_H^h\bm{Q}_H}}{\overline{\bm{Q}}}\odot\frac{(\bm{I}_H^h\bm{Q}_H)\odot\bm{W}_{\bm{z}}}{\bm{Q}}\right\|_{\infty}\,,\\
        &= \left\|\frac{|\bm{I}_H^h\bm{Q}_H|^2}{|\bm{Q}|^2}\odot\frac{|\bm{Q}|^2}{\bm{I}_H^h\bm{I}_h^H|\bm{Q}|^2}\right\|_{\infty}\,,\\
        &= \left\|\frac{\bm{I}_H^h|\bm{Q}_H|^2}{\bm{I}_H^h\bm{I}_h^H|\bm{Q}|^2}\right\|_{\infty}\,,\\
        &= \left\|\bm{I}_H^h\frac{|\bm{Q}_H|^2}{\bm{I}_h^H|\bm{Q}|^2}\right\|_{\infty}\,,\\
        &= \left\|\bm{W}_{\bm{u}}^H\right\|_{\infty}\,.\\ 
    \end{aligned}
\end{equation*}
where the last equality is discussed in the proof of Lemma~\ref{lem:operator_norm}.  Sharpness of the second inequality is straightforward to check by setting both $\bm{Q}$ and $\nabla_{\bm{z}_{k}}\widetilde{\Phi}_k\left(\bm{z}_{k}; \bm{z}_k^{(j)}\right)$ as vectors of ones, in which case $\nabla_{\bm{z}_{H,k}}\widetilde{\Phi}_{H,k}\left(\bm{I}_h^H\left(\bm{W}_{\bm{z}}\odot\bm{z}_{k}\right); \bm{z}_k^{(j)}\right)$ is also a vector of ones.
\end{proof}

\section{Proof for Theorem~\ref{prop:descent_direction}}\label{app:proof_descent_direction}
\begin{proposition}\label{prop:multigrid_convergence}
Suppose that 
\begin{equation*}
    \bm{e}_{H,k} = -\bm{A}\odot\left(\bm{I}_h^H\nabla_{\bm{z}_k} \widetilde{\Phi}_k\left(\bm{z}_k'; \bm{z}_k^{(j)}\right)\right)\,,
\end{equation*} 
for some nonnegative vector $\bm{A}\in\mathbb{R}^{m_H^2}$. Then, $\bm{I}_H^h\left(\bm{e}_{H,k}\right)$ is a descent direction, i.e.
\begin{equation*}                
    \nabla_{\bm{z}_k} \widetilde{\Phi}_k\left(\bm{z}_k'; \bm{z}_k^{(j)}\right)^*\left(\bm{I}_H^h\bm{e}_{H,k}\right)\leq0\,.
\end{equation*}
    Additionally, if $\bm{A}$ is positive, then equality implies $\bm{I}_h^H\nabla_{\bm{z}_k} \widetilde{\Phi}_k\left(\bm{z}_k'; \bm{z}_k^{(j)}\right)=0$.
\begin{proof}
Since $\operatorname{diag}(\bm{A})$ is positive semidefinite, we have
\begin{equation*}
\begin{aligned}
    \nabla_{\bm{z}_k} \widetilde{\Phi}_k\left(\bm{z}_k'; \bm{z}_k^{(j)}\right)^*\left(\bm{I}_H^h\bm{e}_{H,k}\right) &= -\nabla_{\bm{z}_k} \widetilde{\Phi}_k\left(\bm{z}_k'; \bm{z}_k^{(j)}\right)^*\left(\bm{I}_H^h\bm{A}\odot\left(\bm{I}_h^H\nabla_{\bm{z}_k} \widetilde{\Phi}_k\left(\bm{z}_k'; \bm{z}_k^{(j)}\right)\right)\right)\,,\\
    &= -\nabla_{\bm{z}_k} \widetilde{\Phi}_k\left(\bm{z}_k'; \bm{z}_k^{(j)}\right)^*\bm{I}_H^h\operatorname{diag}(\bm{A})\left(\bm{I}_h^H\nabla_{\bm{z}_k} \widetilde{\Phi}_k\left(\bm{z}_k'; \bm{z}_k^{(j)}\right)\right) \,,\\
    &=-4\left(\left(\bm{I}_h^H\nabla_{\bm{z}_k} \widetilde{\Phi}_k\left(\bm{z}_k'; \bm{z}_k^{(j)}\right)\right)\right)^*\operatorname{diag}(\bm{A})\left(\bm{I}_h^H\nabla_{\bm{z}_k} \widetilde{\Phi}_k\left(\bm{z}_k'; \bm{z}_k^{(j)}\right)\right) \,,\\
    &\leq 0\,.
\end{aligned}
\end{equation*}

\end{proof}
\end{proposition}
\subsection{Proof for Theorem~\ref{prop:descent_direction}}\label{app:proof_coarse_grid_update}
    \begin{proof}
Solving~Eqn.~\ref{eqn:coarse_proximal} by setting the complex gradient of the objective function to zero yields:
\begin{equation*}
    \widehat{\bm{z}_{H,k}} = \bm{z}'_{H,k} + \frac{\overline{\bm{Q}_H}}{\bm{u}_H + |\bm{Q}_H|^2}\odot\left(\mathcal{R}_{H,k}\left(\bm{z}_k^{(j)}\right) - \bm{Q}_H\odot\bm{z}'_{H,k}\right)\,,
\end{equation*}
which, by definition, implies
\begin{equation*}
    \bm{e}_{H,k} = -\frac{\overline{\bm{Q}_H}}{\bm{u}_H + |\bm{Q}_H|^2}\odot\left(\bm{Q}_H\odot\bm{z}_{H,k}' - \mathcal{R}_{H,k}\left(\bm{z}_k^{(j)}\right)\right)\,.
\end{equation*}
Then, by Lemma~\ref{lem:commutative} and Lemma~\ref{lem:unary}, we have
\begin{equation*}
        \begin{aligned}
            \bm{e}_{H,k} &= -\frac{\overline{\bm{Q}_H}}{\bm{u}_H + |\bm{Q}_H|^2}\odot\left(\bm{Q}_H\odot\bm{z}_{H,k}' - \bm{I}_h^H\left(\frac{(\bm{I}_H^h\bm{Q}_H )\odot\bm{W}_{\bm{z}}}{\bm{Q}}\odot\mathcal{R}_k\left(\bm{z}_k^{(j)}\right)\right)\right)\,,\\
            &=-\frac{\overline{\bm{Q}_H}}{\bm{u}_H+|\bm{Q}_H|^2}\odot\left(  \bm{Q}_H\odot\bm{I}_h^H(\bm{W}_{\bm{z}}\odot\bm{z}'_{k})-\bm{Q}_H \odot\bm{I}_h^H\left(\frac{\bm{W}_{\bm{z}}}{\bm{Q}}\odot\mathcal{R}_k\left(\bm{z}_k^{(j)}\right)\right)\right)\,,\\
            &=-\frac{|\bm{Q}_H|^2}{\bm{u}_H+|\bm{Q}_H|^2}\odot\bm{I}_h^H\left(\bm{W}_{\bm{z}}\odot\left(\bm{z}'_{k}-\frac{1}{\bm{Q}}\odot\mathcal{R}_k\left(\bm{z}_k^{(j)}\right)\right)\right)\,,\\
            &= -\frac{|\bm{Q}_H|^2}{\bm{u}_H+|\bm{Q}_H|^2}\odot\bm{I}_h^H\left(\frac{\bm{W}_{\bm{z}}}{|\bm{Q}|^2}\odot\nabla_{\bm{z}_k}\widetilde{\Phi}_k\left(\bm{z}_k'; \bm{z}_k^{(j)}\right)\right)\,,\\
            &=-\frac{|\bm{Q}_H|^2}{\bm{u}_H+|\bm{Q}_H|^2}\odot\bm{I}_h^H\left(\frac{1}{\bm{I}_H^h\bm{I}_h^H|\bm{Q}|^2}\odot\nabla_{\bm{z}_k}\widetilde{\Phi}_k\left(\bm{z}_k'; \bm{z}_k^{(j)}\right)\right)
        \,,\\
        &=-\frac{|\bm{Q}_H|^2}{\bm{u}_H+|\bm{Q}_H|^2}\odot\bm{I}_h^H\left(\bm{I}_H^h\left(\frac{1}{\bm{I}_h^H|\bm{Q}|^2}\right)\odot\nabla_{\bm{z}_k}\widetilde{\Phi}_k\left(\bm{z}_k'; \bm{z}_k^{(j)}\right)\right)
        \,,\\
        &=-\frac{|\bm{Q}_H|^2}{\bm{u}_H+|\bm{Q}_H|^2}\odot\frac{1}{\bm{I}_h^H|\bm{Q}|^2}\odot\bm{I}_h^H\left(\nabla_{\bm{z}_k}\widetilde{\Phi}_k\left(\bm{z}_k'; \bm{z}_k^{(j)}\right)\right)
        \,.\\
        \end{aligned}
    \end{equation*}

     The proof follows from  Proposition~\ref{prop:multigrid_convergence} with
    \begin{equation*}
            \bm{A} = \frac{1}{\bm{u}_H+|\bm{Q}_H|^2}\odot\frac{|\bm{Q}_H|^2}{\bm{I}_h^H|\bm{Q}|^2}\,.     
    \end{equation*}
    \end{proof}

\section{Proof for Proposition~\ref{prop:automatic_regularization_selection}}\label{app:proof_automatic_regularization_selection}

\begin{proof}
First, if we solve for the fine-grid update at $\bm{z}_k'$ in the $k$-th scanning region:
\begin{equation*}
    \bm{z}_k^+=\operatornamewithlimits{argmin}_{\bm{z}_k}\widetilde{\Phi}_k\left(\bm{z}_k; \bm{z}_k^{(j)}\right)+\frac{1}{2}\bm{u}^\top\left|\bm{z}_k-\bm{z}_k'\right|^2\,,
\end{equation*}
the resulting update is:
\begin{equation*}
    \bm{e}_k = \bm{z}_k^+-\bm{z}_k' = -\frac{1}{\bm{u}+|\bm{Q}|^2}\odot\nabla_{\bm{z}_k}\widetilde{\Phi}_k\left(\bm{z}_k'; \bm{z}_k^{(j)}\right)\,,
\end{equation*}

The coarse-grid update at $\bm{z}'_{H,k}= \bm I_h^H(\bm{W}_{\bm{z}}\odot\bm{z}'_k)$ (see Eqn.~\ref{eqn:coarse_grid_rPIE_update}) is 
\begin{equation*}
    \bm{e}_{H,k} = -\frac{1}{\bm{u}_H+|\bm{Q}_H|^2}\odot\frac{|\bm{Q}_H|^2}{\bm{I}_h^H|\bm{Q}|^2}\odot\bm{I}_h^H\left(\nabla_{\bm{z}_k}\widetilde{\Phi}_k\left(\bm{z}_k'; \bm{z}_k^{(j)}\right)\right)\,.
\end{equation*}
Hence, both fine-grid and coarse-grid updates are functions of the fine-grid gradient.

Then, with a slight abuse of notation, we define 
$\bm{e}_k(\bm{g})$ and $\bm{e}_{H,k}(\bm{g})$, which map fine-grid gradient $\bm{g}\in\mathbb{C}^{m^2}$ to fine-grid update in $\mathbb{C}^{m^2}$ and coarse-grid update in $\mathbb{C}^{(m/2)^2}$, respectively:
\begin{equation*}
    \bm{e}_k(\bm{g}) = -\frac{1}{\bm{u}+|\bm{Q}|^2}\odot\bm{g}\quad\text{and}\quad\bm{e}_{H,k}(\bm{g})=-\frac{1}{\bm{u}_H+|\bm{Q}_H|^2}\odot\frac{|\bm{Q}_H|^2}{\bm{I}_h^H|\bm{Q}|^2}\odot\bm{I}_h^H\bm{g}
    \,.
\end{equation*}

Since $\|\bm{I}_h^H\|_\infty  = \|\bm{I}_H^h\|_\infty = 1$ by Lemma~\ref{lem:operator_norm}, it is straightforward to see that
\begin{equation*}
    \left\{\bm{I}_h^H\bm{g}:\bm{g}\in\mathbb{C}^{m^2},\|\bm{g}\|_\infty\leq1\right\} = \left\{\bm{g}_H\in\mathbb{C}^{(m/2)^2}:\|\bm{g}_H\|_\infty\leq1\right\}\,.
\end{equation*}
Hence, we observe that for all indices $r$,
\begin{equation*}
\begin{aligned}
    \max_{\|\bm{g}\|_\infty\leq1}\left[\left|\bm{e}_{H,k}(\bm{g})\right|\right]_r &= \max_{\|\bm{g}_H\|_\infty\leq1}\left[\left|\frac{1}{\bm{u}_H+|\bm{Q}_H|^2}\odot\frac{|\bm{Q}_H|^2}{\bm{I}_h^H|\bm{Q}|^2}\odot\bm{g}_H\right|\right]_r\,,\\
    &= \left[\frac{1}{\bm{u}_H+|\bm{Q}_H|^2}\odot\frac{|\bm{Q}_H|^2}{\bm{I}_h^H|\bm{Q}|^2}\right]_r\,,\\
\end{aligned}
\end{equation*}
where the equality for all indices $r$ is achieved when $\bm{g}_H=\bm{J}$, a vector of ones.

Similarly, we have
\begin{equation*}
    \max_{\|\bm{g}\|_\infty\leq1} \left[\left|\bm{I}_h^H\bm{e}_k(\bm{g})\right|\right]_r = \left[\bm{I}_h^H\frac{1}{\bm{u}+|\bm{Q}|^2}\right]_r\,.
\end{equation*}

Finally, by the definition of $\bm{u}_H$ and Jensen's inequality for the convex function $t\mapsto 1/t$ on $(0,\infty)$, applied to each restriction bin, we have:
\begin{equation*}
    \begin{aligned}
        \bm{u}_H &=  \frac{|\bm{Q}_H|^2}{\bm{I}_h^H|\bm{Q}|^2}\odot\bm{I}_h^H\bm{u}\,,\\
        \bm{u}_H + |\bm{Q}_H|^2&= \frac{|\bm{Q}_H|^2}{\bm{I}_h^H|\bm{Q}|^2}\odot\bm{I}_h^H\bm{u}+|\bm{Q}_H|^2\,,\\
        \bm{u}_H + |\bm{Q}_H|^2&= \frac{|\bm{Q}_H|^2}{\bm{I}_h^H|\bm{Q}|^2}\odot\bm{I}_h^H\left(\bm{u}+|\bm{Q}|^2\right)\,,\\
        \bm{I}_h^H\left(\frac{1}{\bm{u}+|\bm{Q}|^2}\right)&\geq \frac{1}{\bm{I}_h^H\left(\bm{u}+|\bm{Q}|^2\right)} = \frac{1}{\bm{u}_H + |\bm{Q}_H|^2}\odot\frac{|\bm{Q}_H|^2}{\bm{I}_h^H|\bm{Q}|^2}\,,\\
        \max_{\|\bm{g}_2\|_\infty\leq1} \left[\left|\bm{I}_h^H\bm{e}_k(\bm{g}_2)\right|\right]_r &\geq \max_{\|\bm{g}_1\|_\infty\leq1}\left[\left|\bm{e}_{H,k}(\bm{g}_1)\right|\right]_r\quad\text{for all }r\,.
    \end{aligned}
\end{equation*}
\end{proof}

\end{document}